\documentclass[]{amsart}
\usepackage{amsmath}
\usepackage{amsfonts}
\usepackage{amssymb}
\usepackage{graphics}
\usepackage{enumerate}

\theoremstyle{plain}
\newtheorem{Thm}{Theorem}[section]

\newtheorem{Lemma}[Thm]{Lemma}
\newtheorem*{MainThm}{The Main Theorem}
\newtheorem{Prop}[Thm]{Proposition}
\theoremstyle{definition}
\newtheorem{Def}[Thm]{Definition}

\theoremstyle{remark}
\newtheorem*{notation}{Notation}
\newtheorem*{remark}{Remark}

\newenvironment{enumeratei}{\begin{enumerate}[\upshape (i)]}{\end{enumerate}}

\def\D{\mathcal{D}}
\def\H{H_1(\Sigma)}
\def\N{\mathbb{N} \cup \{0\}}
\def\R{\mathbb{R}}
\def\Z{\mathbb{Z}}
\def\mtwo{\ (\textnormal{mod}\ 2)}

\begin{document}

\title[Virtual Strings for Closed Curves]{Virtual Strings for Closed Curves with Multiple Components}
\author{William J. Schellhorn}
\address{Department of Mathematics\\
  Louisiana State University\\
  Baton Rouge, LA 70803-4918}
\email{schellho@math.lsu.edu}
\date{March 8, 2005}


\begin{abstract}
A Gauss paragraph is a combinatorial formulation of a generic closed curve with multiple components on some surface.  A virtual string is a collection of circles with arrows that represent the crossings of such a curve.  Every closed curve has an underlying virtual string and every virtual string has an underlying Gauss paragraph.  A word-wise partition is a partition of the alphabet set of a Gauss paragraph that satisfies certain conditions with respect to the Gauss paragraph.  In this paper we use the theory of virtual strings to obtain a combinatorial description of closed curves in the 2-sphere (and therefore $\mathbb{R}^2$) in terms of Gauss paragraphs and word-wise partitions.
\end{abstract}

\maketitle

\tableofcontents

\begin{section}{Introduction}
A \textit{word} in a finite alphabet set, i.e. a finite set with elements referred to as letters, is a sequence of letters up to circular permutation and choice of letters.  A \textit{Gauss paragraph} in a finite alphabet is a finite set of words in the alphabet such that each letter occurs exactly twice in the words of the set.  A Gauss paragraph with only one word is called a \textit{Gauss word} or \textit{Gauss code}.  Such words were introduced by C.F. Gauss as a combinatorial formulation of closed curves on $\R^2$. 

The concept of a virtual string was first introduced by Vladimir Turaev in \cite{Tvs} to describe a single copy of $S^1$ with distinguished ordered pairs of points, which can be represented as a set of arrows attached to the circle.  We extend this definition to allow multiple copies of $S^1$.  Every virtual string has an underlying Gauss paragraph, with the circles and arrows of the virtual string corresponding to the words and letters of the Gauss paragraph, respectively.

A \textit{closed curve} with $N \in \mathbb{N}$ components on a surface $\Sigma$ is a generic smooth immersion of $N$ oriented circles into the surface $\Sigma$.  Every closed curve $\rho$ has an underlying virtual string, with the components and crossings of the closed curve corresponding to the circles and arrows of the virtual string, respectively.  The Gauss paragraph of a closed curve is the Gauss paragraph of its underlying virtual string.

Let $\rho$ be a closed curve on a surface, with crossings labeled by the elements of some finite set $E$.  Pick a base point on each component of $\rho$ that is not a crossing of $\rho$.  Having selected these base points, it is relatively easy to determine the Gauss paragraph $p$ of $\rho$ in the alphabet $E$.  The word of $p$ corresponding to a given component is the sequence of letters formed in the following way: start at the base point of the component and record, in order, the labels of the crossings encountered while traversing the component according to its orientation.  A Gauss paragraph $p$ is said to be \textit{realizable} by a closed curve on a surface $\Sigma$ if there exists some closed curve $\rho$ on the surface $\Sigma$ whose Gauss paragraph is $p$.

Many mathematicians have questioned when a given Gauss word is realizable by a closed curve on $\R^2$, or equivalently on the 2-sphere $S^2$.  Four of the most notable studies include \cite{Gw}, \cite{LMafs}, \cite{Rsad}, and \cite{DTcok}.  In this paper we address when a Gauss paragraph is realizable by a closed curve on $S^2$.  We use the theory of virtual strings to obtain a combinatorial description of closed curves on the 2-sphere in terms of Gauss paragraphs and ``word-wise partitions".

This paper is organized as follows.  In Section 2, we discuss Gauss paragraphs and define the concept of a word-wise partition.  Section 3 contains the definitions for several maps that are related to Gauss paragraphs.  We define virtual strings on multiple copies of $S^1$ in Section 4, and explain how virtual strings have underlying Gauss paragraphs.  In Section 5, we discuss closed curves on surfaces, explain how closed curves have underlying virtual strings, and describe how to construct a surface containing a closed curve realizing a given virtual string.  In Section 6, we observe that a given virtual string may or may not give rise to a word-wise partition, and explain why it is acceptable to restrict our attention to those virtual strings that do give rise to word-wise partitions.  We discuss a homological intersection pairing in Section 7 and relate the pairing to the maps defined in Section 3.  Section 8 contains a statement of the main theorem and the proofs of several lemmas.  These lemmas are used to prove the main theorem in Section 9.  Section 10 contains some additional results about Gauss paragraphs.
\end{section}

\begin{section}{Gauss Paragraphs}
Let us impose the condition that it is not possible to partition the words of a Gauss paragraph into two sets such that the words in one set have no letters in common with the words of the other.  This condition assures that a Gauss paragraph is not a ``disjoint union" of other paragraphs.  A letter in a word $w$ of a Gauss paragraph is called a \textit{single letter} of $w$ if it occurs once in $w$, and a \textit{double letter} of $w$ if it occurs twice in $w$.  Given a set $A$, denote its cardinality by $\#(A)$.  Then the \textit{length} of a word $w$ of a Gauss paragraph is defined as $$2(\#\{\textnormal{double letters in}\ w\}) + \#\{\textnormal{single letters in}\ w\}.$$

Although a Gauss paragraph is by definition a finite set of words, we will write a Gauss paragraph as a finite sequence of words when it is necessary to reference the words individually.  This convention is reasonable because we will not discuss any notion of equivalence amongst Gauss paragraphs.

Throughout the remainder of this section, let $p$ be a Gauss paragraph in an alphabet $E$.

\begin{Def}
A \textit{word-wise partition} $P$ of $E$ with respect to $p$ is a partition of the letters of $E$ that satisfies the following conditions:
\begin{enumeratei}
    \item{$P$ associates two disjoint, possibly empty subsets $A$ and $A'$ to every word $w$ in $p$;}
    \item{for each word $w$ in $p$, the union $$(A \cap \{\textnormal{double letters of}\ w\}) \cup (A' \cap \{\textnormal{double letters of}\ w\})$$ gives a bipartition of the set of double letters of $w$; and}
    \item{if a word $w$ of $p$ has $2n$ letters in common with another word $w'$ of $p$, then $n$ of these letters appear in the two sets in $P$ associated to $w$ and $n$ of them appear in the two sets in $P$ associated to $w'$.}
\end{enumeratei}
\end{Def}
In practice, if $p$ is written as a sequence $(v_1,v_2,...,v_N)$ of words, then we will express a word-wise partition $P$ using the notation $(A_1 \cup A'_1,A_2 \cup A'_2,...,A_N \cup A'_N)$, where the sets $A_n$ and $A'_n$ are associated to the word $v_n$.  Notice that a word-wise partition is a partition of $E$ such that its sets satisfy the above conditions with respect to the words of $p$.  However, it does not partition the letters in each word of $p$ unless $p$ has only one word, and in this case it is a bipartition of the letters in the word (it is the bipartition Turaev defined in \cite{Tvs}).

When we use the terms ``sequence" and ``subsequence", we mean finite sequences of letters that are not considered up to circular permutation.  Recall, however, that words are considered up to circular permutation.  Given a word $w$, if we write $w = i \cdots j$ for some sequence of letters $i \cdots j$, we mean that $w$ can be written in the form $i \cdots j$ up to circular permutation.  A finite sequence $x_1$ of letters is called a $subsequence$ of a word $w$ if $w = x_1x_2$ for some sequence $x_2$ of letters.  Consequently, $x_2$ is also a subsequence of $w$.  The \textit{length} of a sequence is the number of entries in the sequence, and the length of a sequence $s$ will be denoted by $\ell(s)$.  We will use the symbol $\emptyset$ to denote an empty sequence and define $\ell(\emptyset) = 0$.

\begin{Def}
Let $i \in E$ be a double letter of a word $w$ of $p$.  The two \textit{p-sets} of $i$, denoted $p_i$ and $p'_i$, are the two sets of letters in $w$ defined as follows.  Since $w$ is considered up to circular permutation, assume it has the the form $i x_1 i x_2$, with $x_1$ and $x_2$ being subsequences (possibly empty) of the word $w$.  The two $p$-sets are the set of letters occuring exactly once in $x_1$ and the set of letters occuring exactly once in $x_2$.
\end{Def}
Note that $p_i = p'_i$ if $p$ has only one word.  Turaev denotes this set by $w_i$ in \cite{Tvs}, and uses the concept of interlacing to define $w_i$.  In general, two distinct letters $i,j \in E$ in a word $w$ of a Gauss paragraph are called \textit{w-interlaced} if $w$ has the form $i \cdots j \cdots i \cdots j \cdots$ up to circular permutation.  Then $w_i$ is the set of letters $w$-interlaced with $i$.

\begin{Def}
A word-wise partition $P$ of $E$ (with respect to $p$) is \textit{compatible} with $p$ if the following two conditions are satisfied.  Suppose $i,j \in E$ are $w$-interlaced in a word $w$ of $p$.  Then $w=ix_1jx_2ix_3jx_4$ for some (possibly empty) subsequences $x_1,x_2,x_3,x_4$ of $w$.  The first condition is that for any two such letters $i$ and $j$,
$$\begin{aligned}
&\#(w_i \cap w_j) + \#\{\textnormal{single letters of $w$ in}\ x_1\} \\
&\equiv \#(w_i \cap w_j) + \#\{\textnormal{single letters of $w$ in}\ x_3\} \mtwo \\
&\equiv
\begin{cases}
0 \mtwo \hspace{.1in} \textnormal{if $i,j$ appear in different subsets of $P$;} \\
1 \mtwo \hspace{.1in} \textnormal{if $i,j$ appear in the same subset of $P$.}
\end{cases}
\end{aligned}$$
Suppose $i,j \in E$ are single letters in a word $w$ of $p$ that appear in the union of the two subsets associated to $w$ in P.  Then $w=ix_1jx_2$ for some (possibly empty) subsequences $x_1,x_2$ of $w$.  The second condition is that for any two such letters $i$ and $j$,
$$\begin{aligned}
\ell(x_1) &\equiv \ell(x_2) \mtwo\\
&\equiv
\begin{cases}
0 \mtwo \hspace{.1in} \textnormal{if $i,j$ appear in different subsets of $P$;} \\
1 \mtwo \hspace{.1in} \textnormal{if $i,j$ appear in the same subset of $P$.}
\end{cases}
\end{aligned}$$
\end{Def}
\end{section}

\begin{section}{Maps Related to Gauss Paragraphs}
Throughout this section, let $p = (v_1,v_2,...,v_N)$ be a Gauss paragraph in an alphabet set $E$ with word-wise partition $P = (A_1 \cup A'_1,A_2 \cup A'_2,...,A_N \cup A'_N)$.
\begin{notation}
Given a sequence of letters $i \cdots j$, let $o(i \cdots j)$ denote the letters that occur exactly once between $i$ and $j$.  Given a word $w$ of a Gauss paragraph, let $o(w)$ denote the single letters in $w$.
\end{notation}
\begin{Def}
Let $a \cdots b$ and $y \cdots z$ be subsequences of a word $v_n$ of $p$, where $y$ and $z$ are distinct single letters of $v_n$.  Set
$$\delta_n(a \cdots b,\emptyset) = \delta_n(\emptyset,y \cdots z) = 0.$$
Define $\delta_n(a \cdots b,y \cdots z)$ as follows:
\begin{enumeratei}
\item{if $v_n = a \cdots y \cdots b \cdots z \cdots$, set it equal to
$$\#(o(a \cdots y) \cap o(y \cdots b)) + \#(o(a \cdots y) \cap o(b \cdots z)) + \#(o(y \cdots b) \cap o(b \cdots z))$$}
\item{if $v_n = a \cdots z \cdots b \cdots y \cdots$, set it equal to
$$\#(o(a \cdots z) \cap o(y \cdots a)) + \#(o(z \cdots b) \cap o(y \cdots a)) + \#(o(z \cdots b) \cap o(a \cdots z))$$}
\item{if $v_n = a \cdots y \cdots z \cdots b \cdots$, set it equal to
$$\#(o(a \cdots y) \cap o(y \cdots z)) + \#(o(z \cdots b) \cap o(y \cdots z))$$}
\item{if $v_n = a \cdots z \cdots y \cdots b \cdots$, set it equal to
$$\begin{aligned}
\#(o(a \cdots &z) \cap o(b \cdots a)) + \#(o(z \cdots y) \cap o(y \cdots b)) + \#(o(z \cdots y) \cap o(b \cdots a)) \\
&+ \#(o(z \cdots y) \cap o(a \cdots z)) + \#(o(y \cdots b) \cap o(b \cdots a))
\end{aligned}$$}
\item{if $v_n = a \cdots b \cdots y \cdots z \cdots$, set it equal to
$$\#(o(a \cdots b) \cap o(y\cdots z))$$}
\item{if $v_n = a \cdots b \cdots z \cdots y \cdots$, set it equal to
$$\#(o(a \cdots b) \cap o(y \cdots a)) + \#(o(a \cdots b) \cap o(b \cdots z))$$}
\end{enumeratei}
Although $y \ne z$ by hypothesis, no other pairwise distinct requirements are enforced in this definition.  For example, if $a=y$ and $b=z$, then the first, second, third, and sixth cases are the same.  The definition of $\delta_n$ is consistent under all such circumstances.
\end{Def}

\begin{Def}
Let $a x_1 b$ and $y x_2 z$ be two subsequences of $v_n$, where $y$ and $z$ are distinct single letters of $v_n$ and $x_1,x_2$ are subsequences of $v_n$.  Set 
$$\epsilon_n^P(a x_1 b,\emptyset) = \epsilon_n^P(\emptyset,y x_2 z) = 0.$$
If $a=b$ is a double letter of $v_n$, define $\epsilon_n^P(a x_1 b,y x_2 z)$ to be 1 if an odd number of the following statements is true and 0 otherwise:
\begin{enumeratei}
\item{$a=b$ occurs exactly once in $x_2$}
\item{$y$ occurs in $x_1$ and $y \notin A_n \cup A'_n$}
\item{$z$ occurs in $x_1$ and $z \in A_n \cup A'_n$}
\end{enumeratei}
If $a$ and $b$ are distinct single letters of $v_n$, define $\epsilon_n^P(a x_1 b,y x_2 z)$ to be 1 if an odd number of the following statements is true and 0 otherwise:
\begin{enumeratei}
\item{$a$ occurs in $x_2$ and $a \in A_n \cup A'_n$}
\item{$b$ occurs in $x_2$ and $b \notin A_n \cup A'_n$}
\item{$y$ occurs in $x_1$ and $y \notin A_n \cup A'_n$}
\item{$z$ occurs in $x_1$ and $z \in A_n \cup A'_n$}
\end{enumeratei}
Although $y \ne z$ by hypothesis, no other pairwise distinct requirements are enforced in this definition.
\end{Def}
\begin{Def}
Let $X = (x_1,x_2,...,x_N)$ be a sequence consisting of subsequences of the words in $p$, where $x_n$ is a subsequence of the word $v_n$ and more than one $x_n$ is nonempty.  Let $X'$ be the set containing the entries of $X$ that are nonempty subsequences, and define $M = \#(X') \geq 2$.  The set $X$ is called a \textit{cyclic sequence} associated to $p$ if a sequence $(x'_1,x'_2,...,x'_M)$ including all elements in $X'$ can be constructed such that the first letter of $x'_1$ is the last letter of $x'_M$ and the the first letter of $x'_{m}$ is the last letter of $x'_{m-1}$ for $2 \leq m \leq M$.  Note that the first and last letters of each subsequence in a cyclic sequence are single letters of the same word.  Let $\D_p$ denote the collection of all cyclic sequences associated to $p$.
\end{Def}
Suppose $d = (x_1,x_2,...,x_N) \in \D_p$.  Denote the sequence $x_n \in d$ by $d(n)$.  Then $o(d(n))$ is the set of letters in $d(n) = x_n$ that occur exactly once between the first and last letters if it is a nonempty sequence, and the empty set otherwise.  We will now define a collection of maps that involve $\D_p$.  These maps will eventually be related to a homological intersection form.
\begin{Def}
Suppose $i,j$ are distinct single letters in a word $v_n$ of $p$.  Define $\gamma_n^P(i,j)$ to be 0 if exactly one of the letters $i,j$ appears in $A_n \cup A'_n$, and 1 if both or neither of these letters appears in $A_n \cup A'_n$.
\end{Def}
\begin{Def}
The map $W: p \times \D_p \to \N$ is defined by 
$$W(v_n,d) = \sum_{k \neq n} \#(o(v_n) \cap o(d(k)))$$
if $d(n)$ is an empty sequence, and
$$W(v_n,d) = \#(o(ix_1j) \cap o(jx_2i)) + \gamma_n^P(i,j) + \sum_{k \neq n} \# (o(v_n) \cap o(d(k)))$$
if $d(n)$ is subsequence $ix_1j$ of $v_n = i x_1 j x_2$ where $i,j$ are distinct single letters of $v_n$ and $x_1,x_2$ are subsequences of $v_n$.
\end{Def}
\begin{Def}
Suppose $i \in E$ is a double letter in word $v_n$ of $p$.  Then $v_n$ can be written as $ix_1ix_2$ and $ix_2ix_1$ for some subsequences $x_1$ and $x_2$, but we may assume the $p$-sets $p_i$ and $p'_i$ contain the letters that occur exactly once in $x_1$ and $x_2$, respectively.  Then for $d \in \D_p$, set
$$Q_n(p_i,d) = \delta_n(ix_1i,d(n)) + \epsilon_n^P(ix_1i,d(n)) + \sum_{k \neq n}\#(p_i \cap o(d(k)))$$
and 
$$Q_n(p'_i,d) = \delta_n(ix_2i,d(n)) + \epsilon_n^P(ix_2i,d(n)) + \sum_{k \neq n}\#(p'_i \cap o(d(k))).$$
\end{Def}
\begin{Def}
Define a map $D_n : \D_p \times \D_p \to \N$ as follows.  For $d_1, d_2 \in \D_p$, define 
$$D_n(d_1,d_2) = \delta_n(d_1(n),d_2(n)) + \epsilon_n^P(d_1(n),d_2(n)) + \sum_{k \neq n}\#(o(d_1(n)) \cap o(d_2(k))).$$
\end{Def}
\begin{Def}
The word-wise partition $P$ is $compatible$ with $\D_p$ if the following conditions are satisfied:
\begin{enumeratei}
	\item{for every word $v_n$ of $p$, $W(v_n,d) \equiv 0 \mtwo$ for all $d \in \D_p$;}
	\item{if $i \in E$ is a double letter of word $v_n$ of $p$, then 
		$Q_n(p_i,d) \equiv Q_n(p'_i,d) \equiv 0 \mtwo$ for all $d \in \D_p$; and}
	\item{$\sum_{n=1}^N D_n(d_1,d_2) \equiv 0 \mtwo$ for all $d_1,d_2 \in \D_p$.}
\end{enumeratei}
\end{Def}
\end{section}

\begin{section}{Virtual Strings}
Turaev introduced the term ``virtual string" in \cite{Tvs} to describe a single copy of $S^1$ with distinguished ordered pairs of points, which can be represented as a set of arrows attached to the circle.  We extend this definition to allow multiple copies of $S^1$.  Consider a collection of $N$ copies of $S^1$ with a distinguished set of points on each copy, such that the total number of distinguished points amongst all $N$ copies is even.  A \textit{virtual string} $\alpha$ with $N$ components is such a collection with the distinguished points partitioned into ordered pairs.  The copies of $S^1$ in $\alpha$ are called the \textit{core circles} of $\alpha$ and the ordered pairs are called the $arrows$ of $\alpha$.  For an arrow $(a,b)$ of $\alpha$, the endpoints $a$ and $b$ are called its $tail$ and $head$, respectively.  We will impose the condition that it is not possible to partition the core circles of a virtual string into two sets such that the circles in one set have no arrows in common with the circles of the other.  This condition assures that a virtual string is not a ``disjoint union" of other virtual strings.  Figure \ref{F:VSExample} depicts a virtual string with three core circles and twelve arrows.

\begin{figure}
\scalebox{.28}{\includegraphics{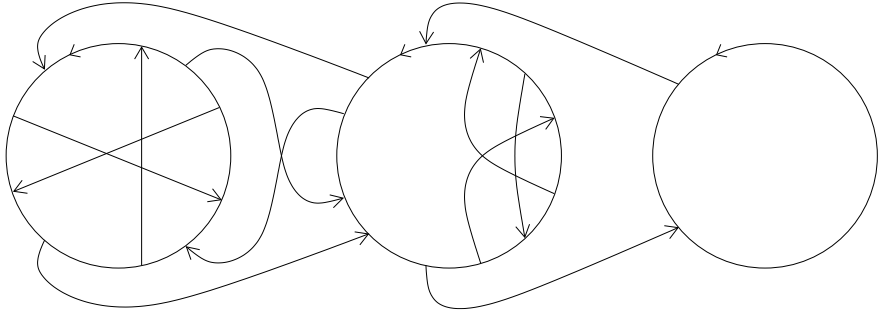}}
\caption{A virtual string with three core circles and twelve arrows.}
\label{F:VSExample}
\end{figure}

Let $S$ be a core circle of a virtual string $\alpha$.  Two distinguished points $a,b$ on $S$ separate $S$ into two arcs, namely $ab$ and $ba$.  Suppose $e = (a,b)$ is an arrow of $\alpha$.  A different arrow $f = (c,d)$ on $S$ is said to \textit{link} $e$ if one endpoint of $f$ lies on the interior $(ab)^\circ$ of the arc $ab$ and the other lies on the interior $(ba)^\circ$ of the arc $ba$.  The arrow $f$ links $e$ \textit{positively} if the arrow endpoints lie in the cyclic order $a,d,b,c$ around $S$, and the arrow $f$ links $e$ \textit{negatively} if the arrow endpoints lie in the cyclic order $a,c,b,d$ around $S$.  Note that if $f$ links $e$ positively then $e$ links $f$ negatively.

Two virtual strings are $homeomorphic$ if there is an orientation-preserving homeomorphism of their core circles such that the arrows of the first string are mapped onto the arrows of the second string.  The homeomorphism classes of virtual strings will also be called virtual strings.

Every virtual string $\alpha$ has an underlying Gauss paragraph $p = p_\alpha$.  Pick any alphabet set $E$ with $\#(E)$ equal to the number of arrows in $\alpha$, and label each arrow with a different letter in $E$.  Select a base point on each core circle of $\alpha$ that is not an endpoint of any arrow.  The word of $p$ associated to the core circle $S$ of $\alpha$ is obtained as follows: starting at the base point of $S$, traverse $S$ in the positive direction and record the label of an arrow each time one of its endpoints is encountered.  The resulting word will be well-defined up to circular permutations (and of course the choice of arrow labels).  The Gauss paragraph underlying $\alpha$ is the set of words obtained in this way, with one word associated to every core circle.
\end{section}

\begin{section}{Closed Curves on Surfaces}
Recall that a \textit{surface} is a smooth oriented 2-dimensional manifold.  Also recall that a smooth map from a collection of oriented circles $\amalg S^1$ into a surface $\Sigma$ is called an \textit{immersion} if its differential is nonzero at all points of the circle.  For an immersion $\rho : \amalg S^1 \to \Sigma$, a point $x \in \Sigma$ with $\#(\rho^{-1}(x))=2$ is called a \textit{double point} or \textit{crossing} of $\rho$.  The immersion $\rho$ is called \textit{generic} if $\#(\rho^{-1}(x)) \leq 2$ for all $x \in \Sigma$, it has a finite set of double points, and all its double points are transverse intersections of two branches.  A generic smooth immersion of $N$ oriented circles into a surface $\Sigma$ is called a \textit{closed curve} with $N$ components on the surface $\Sigma$.

Every closed curve $\rho: \amalg S^1 \to \Sigma$ has an underlying virtual string $\alpha = \alpha_\rho$.  The core circles of $\alpha$ are the copies of $S^1$ in the domain of $\rho$.  The arrows of $\alpha$ are all ordered pairs $(a,b)$ of distinguished points such that $\rho(a) = \rho(b)$ and the pair (a positive tangent vector of $\rho$ at $a$, a positive tangent vector of $\rho$ at $b$) is a positive basis in the tangent space of $\rho(a)$.  A virtual string is said to be \textit{realized} by a closed curve $\rho : \amalg S^1 \to \Sigma$ if it is homeomorphic to $\alpha_\rho$.  The virtual string in Figure \ref{F:VSExample} is realized by the closed curve in Figure \ref{F:CCExample}, which can be considered as a closed curve in either $\R^2$ or $S^2$.  The Gauss paragraph of a closed curve is the Gauss paragraph of its underlying virtual string.

\begin{figure}
\scalebox{.38}{\includegraphics{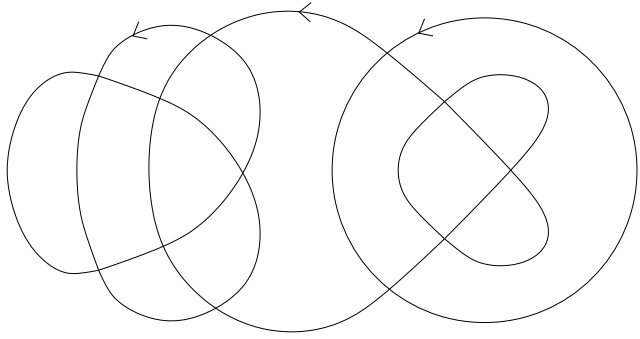}}
\caption{A closed curve with three components and twelve crossings that realizes the virtual string in Figure \ref{F:VSExample}.}\label{F:CCExample}
\end{figure}

Every virtual string admits a canonical realization by a closed curve on a surface.  Turaev showed this result for virtual strings with one core circle in \cite{Tvs}.  He used a well-known construction of surfaces from four-valent graphs, see for example \cite{Ccic}.  We will describe the same construction for virtual strings with multiple core circles.  Let $\alpha$ be a virtual string with more than one core circle.  Transform $\alpha$ into a 1-dimensional CW-complex $\Gamma = \Gamma_\alpha$ by identifying the head and tail of each arrow in $\alpha$.  A thickening of $\Gamma$ gives a surface $\Sigma_\alpha$ in the following manner.  The 0-cells of $\Gamma$ are 4-valent vertices.  A vertex $v \in \Gamma$ results from an arrow $(a,b)$, with points $a$ and $b$ on some core circles of $\alpha$.  Call these core circles $S_a$ and $S_b$, and note that possibly $S_a=S_b$.  A neighborhood of a point $x$ in a core circle is an oriented arc, which $x$ splits into one incoming arc and one outgoing arc with respect to $x$.  The incoming and outgoing arcs in neighborhoods of $a$ in $S_a$ and $b$ in $S_b$ can be identified with the four arcs in a neighborhood of $v$ in $\Gamma$.  Therefore this neighborhood of $v$ can be embedded in the unit 2-disc $D^2 = \{(x,y) \in \R^2 |\ x^2 + y^2 \leq 1\}$ so that the image of $v$ is the origin, the images of the incoming arcs at $a$ and $b$ are the intervals $0 \times [-1,0]$ and $[-1,0] \times 0$ respectively, and the images of the outgoing arcs at $a$ and $b$ are the intervals $0 \times [0,1]$ and $[0,1] \times 0$ respectively.  By repeating this procedure for all vertices of $\Gamma$, the vertices can be thickened to 2-discs endowed with counterclockwise orientation.  A 1-cell of $\Gamma$ either connects two different vertices or forms a loop at a single vertex.  The 1-cells of $\Gamma$ can be thickened to ribbons, with the thickening uniquely determined by the condition that the orientations of the 2-discs extend to their unions with the ribbons.  Thickening $\Gamma$ in the way just described gives an embedding of $\Gamma$ onto a surface $\Sigma_\alpha$.  Notice that $\Sigma_\alpha$ is a compact, connected, oriented surface with boundary.  A closed curve $\rho_\alpha : \amalg S^1 \to \Sigma_\alpha$ realizing $\alpha$ is obtained by composing the natural projection $\amalg S^1 \to \Gamma$ with the inclusion $\Gamma \hookrightarrow \Sigma_\alpha$.

The surface $\Sigma_\alpha$ constructed above is the surface of minimal genus containing a closed curve realizing $\alpha$.  Suppose $\rho : \amalg S^1 \hookrightarrow \Sigma$ is a generic closed curve realizing $\alpha$ on some surface $\Sigma$.  Then a regular neighborhood of $\rho(\amalg S^1)$ in $\Sigma$ is homeomorphic to $\Sigma_\alpha$.  This homeomorphism can be chosen to transform $\rho$ into $\rho_\alpha$, that is composing $\rho_\alpha$ with an orientation-preserving embedding of $\Sigma_\alpha$ into $\Sigma$ results in the curve $\rho$.  Gluing 2-discs to all the components of the boundary of $\Sigma_\alpha$ gives a closed surface of minimal genus that contains a curve realizing $\alpha$.

\begin{notation}
Let $\rho : \coprod_{n=1}^N S'_n \to \Sigma$ be a closed curve with $N \in \mathbb{N}$ components on some surface $\Sigma$, where each $S'_n$ is a copy of $S^1$.  Then $\rho$ has an underlying virtual string $\alpha$ with core circles $S_1,S_2,...,S_N$, where the core circle $S_n$ corresponds to the circle $S'_n$ in the domain of $\rho$.  Let $p = (v_1,v_2,...,v_N)$ be the underlying Gauss paragraph of $\alpha$, with word $v_n$ corresponding to core circle $S_n$.  A cyclic sequence $d = (x_1,x_2,...,x_N) \in \D_p$ corresponds to a sequence $C_d = (c_1,c_2,...,c_N)$ of arcs in $\alpha$ such that the arc $c_n \subset S_n$ corresponds to the sequence $x_n$ in $d$.  Moreover, $C_d$ corresponds to a sequence $C'_d = (c'_1,c'_2,...,c'_N)$ of arcs on the components $S'_n$ in the domain of $\rho$ such that the arc $c'_n \subset S'_n$ corresponds to the arc $c_n \subset S_n$.  Let $C_d(n)$ denote the arc $c_n$ of $C_d$ and $C'_d(n)$ denote the arc $c'_n$ of $C'_d$.  In addition, let $\rho(C'_d)$ denote the loop on $\Sigma$ that is the union $\bigcup_{n=1}^N \rho(C'_d(n))$.  Table \ref{Ta:Correspondences} displays the correspondences and notation discussed thus far.
\end{notation}

\begin{table}
\caption{Summary of notational conventions}
\begin{tabular}{ccc}
\hline
Gauss paragraph $p$  &Virtual string $\alpha$  &Closed curve $\rho$  \\
\hline
word $v_n$  &core circle $S_n$  &component $\rho(S'_n)$  \\
\hline
letter $i$ in $v_n$  &arrow $e_i = (a_i,b_i)$ &crossing $i$ at $\rho(a_i)=\rho(b_i)$ \\
\hline
$p$-sets $p_i$ and $p'_i$  &arcs $a_ib_i$ and $b_ia_i$  
& two loops in $\rho(S'_n)$ based at $\rho(a_i)=\rho(b_i)$  \\
\hline
sequence $d \in \D_p$  &sequence of arcs $C_d$
&loop $\rho(C'_d)$ with segments $\rho(C'_d(n))$ \\
\hline
\end{tabular}
\label{Ta:Correspondences}
\end{table}
\end{section}

\begin{section}{Observations and Restrictions}
Every virtual string with an even number of arrow endpoints on each of its core circles naturally gives rise to a word-wise partition.  Let $\alpha$ be such a virtual string, and suppose $S$ is a core circle of $\alpha$.  Denote the set of arrows of $\alpha$ with tails on $S$ by $arr_t(S)$.  For $e = (a,b)$ and $f = (c,d)$ in $arr_t(S)$, define $q(e,f) \in \Z$ to be the number of arrowheads lying on the semi-open arc $ac - \{a\} \subset S$ minus the number of arrowtails lying on $ac - \{a\}$.  Set $q(e,f) = 0$ when $e = f$.  Use $q$ to define an equivalence relation on $arr_t(S)$ by defining two arrows $e = (a,b)$ and $f = (c,d)$ to be \textit{equivalent} if $q(e,f) \equiv 0 \mtwo$.  Then the arrows $e$ and $f$ are equivalent if either $e = f$ or the number of arrow endpoints lying on the interior of the arc $ac \subset S$ is odd.  Notice that there are at most two equivalence classes on $arr_t(S)$ that result from this relation, and these classes partition $arr_t(S)$ into two subsets.  Pick any alphabet set $E$ with $\#(E)$ equal to the number of arrows in $\alpha$, and label each arrow with a different letter in $E$.  Suppose $E' \subset E$ contains the letters that label the arrows in $arr_t(S)$.  Then the bipartition of $arr_t(S)$ induces a bipartition of the set $E'$.  Applying the same procedure on all core circles of $\alpha$ creates a word-wise partition of $E$.  Thus $\alpha$ naturally gives rise to a pair (the underlying Gauss word $p$ of $\alpha$, a word-wise partition of $E$ with respect to $p$). 

If $\alpha$ is the underlying virtual string of a closed curve $\rho$ on some surface, then the above construction gives a partition of the set of labels (namely $E$) on the double points of $\rho$.  Thus $\rho$ gives rise to the pair (the underlying Gauss word $p$ of $\alpha$, a word-wise partition of $E$ with respect to $p$) as well.  Note that each component of $\rho$ has an even number of intersections between it and the other components because we assumed each core circle of $\alpha$ had an even number of arrow endpoints.

An important observation should be made here.  Let $e$ be an arrow of $\alpha$ with tail on one core circle, say $S_1$, and head on a different core circle, say $S_2$.  Then $e$ corresponds to a letter $i \in E$ that appears in one of the sets associated to $S_1$ in the word-wise partition, and in neither of the sets associated to $S_2$.  Moreover, the sign of crossing $i$ of $\rho$ is +1 with respect to the component of $\rho$ corresponding to $S_2$ and -1 with respect to the component of $\rho$ corresponding to $S_1$.

A pair (a Gauss paragraph $p$ in an alphabet $E$, a word-wise partition $P$ of $E$) is said to be \textit{realizable} by a closed curve on a surface if there exists a closed curve on the surface that gives rise to the pair.  The surface of particular interest to us is $S^2$, but observe that a pair (a Gauss paragraph $p$ in an alphabet $E$, a word-wise partition $P$ of $E$) is realizable by a closed curve on $S^2$ if and only if the pair is realizable by a closed curve on $\R^2$.

If a virtual string has a core circle with an odd number of arrow endpoints on it, then the pairing $q(e,f)$ above is not well-defined on that core circle.  Therefore such virtual strings do not give rise to word-wise partitions in any sense that we will discuss in this paper.  It is well known that a closed curve on $S^2$ cannot have an odd number of intersections between any two distinct components, and therefore each component of a closed curve on $S^2$ cannot have an odd number of intersections between it and all the other components.  Consequently, the underlying virtual string of a closed curve on $S^2$ cannot have a core circle with an odd number of arrow endpoints on it.

Since Gauss paragraphs that can be realized by closed curves on $S^2$ are our primary interest, it is reasonable to now restrict our attention to virtual strings with an even number of arrow endpoints on each core circle and closed curves for which each component has an even number of crossings between it and all the other components.
\end{section}

\begin{section}{Homological Intersection Form}
In this section, we relate a homological intersection pairing to the maps described in Section 3.  We develop formulas that help prove the main theorem.

Let $\Sigma$ be an oriented surface, and let $\H$ denote its first integral homology group $H_1(\Sigma;\Z)$.  Let $\rho : \coprod_{n=1}^N S'_n \to \Sigma$ be a closed curve with $N \in \mathbb{N}$ components, where each $S'_n$ is a copy of $S^1$.  Then $\rho$ has an underlying virtual string $\alpha$, so suppose the Gauss paragraph underlying $\alpha$ is $p = (v_1,v_2,...,v_N)$.  Label the $N$ core circles of $\alpha$ by $S_1,S_2,...,S_N$, with core circle $S_n$ corresponding to circle $S'_n$ and word $v_n$.  A distinguished point $z$ on $S_n$ corresponds to a point on the circle $S'_n$, and for convenience this point on $S'_n$ will also be called $z$.
 
If $e = (a,b)$ is an arrow on core circle $S_n$ of $\alpha$, then $\rho(a)=\rho(b)$.  Therefore the images $\rho(ab)$ and $\rho(ba)$ of the arcs $ab$ and $ba$ of $S'_n$, respectively, are loops in $\Sigma$.  To stay consistent with Turaev's convention in \cite{Tvs}, we will use the notation $[e]$ and $[e]^*$ to denote the homology classes $[\rho(ab)]$ and $[\rho(ba)]$ of $\H$, respectively.  The image under $\rho$ of each circle $S'_n$ is also a loop in $\Sigma$, so let $[S_n]$ denote the homology class $[\rho(S'_n)] \in \H$.  Let $(x'_1,x'_2,...,x'_M)$ with $M \geq 2$ be a sequence of arcs on the circles $S'_n$ that satisfies the following properties: at most one arc exists on each $S'_n$, the initial point of $x'_1$ and the terminal point of $x'_M$ have the same image under $\rho$, and the initial point of $x'_m$ and the terminal point of $x'_{m-1}$ have the same image under $\rho$ for $2 \leq m \leq M$.  Then the union $\bigcup_{m=1}^M \rho(x'_m)$ constitutes a loop in $\Sigma$.  The sequence of arcs $(x'_1,x'_2,...,x'_M)$ corresponds to sequence of arcs on the core circles $S_n$, which in turn corresponds to a sequence $d = (x_1,x_2,...,x_N)$ with $x_n$ being a subsequence (possibly empty) of the word $v_n$ of $p$.  Since more than one $x_n$ is nonempty, it follows that $d \in \D_p$.  Therefore let $[C_d]$ denote the homology class of the union $\bigcup_{m=1}^M \rho(x'_m) = \bigcup_{n=1}^N \rho(C'_d(n)) = \rho(C'_d)$ in $\H$.  Notice that if a loop in $\rho$ has lone segments on multiple components and the orientations of these segments agree, then the loop can be expressed as such a union for some cyclic sequence $d \in \D_p$.

In accordance with the restrictions discussed at the end of Section 6, assume that each component of $\rho$ has an even number of crossings between it and all the other components.  Then $\alpha$ has an even number of arrow endpoints on each core circle.  To this point, we have not specified an alphabet set for $p$.  Suppose $p$ is a paragraph in the alphabet set $E$, and suppose $\alpha$ gives rise to the word-wise partition $P = (A_1 \cup A'_1, A_2 \cup A'_2,..., A_N \cup A'_N)$ of $E$.  Let $e_i = (a_i,b_i)$ denote the arrow in virtual string $\alpha$ that corresponds to the letter $i \in E$.  Label a crossing of the closed curve $\rho$ with the letter $i$ if the crossing corresponds to the arrow $e_i$.  Under these notational conventions, we may assume that $p_i$ and $p'_i$ contain the letters corresponding to the arrows in $\alpha$ with exactly one endpoint on $(a_ib_i)^\circ$ and $(b_ia_i)^\circ$, respectively.

The orientation of the surface $\Sigma$ determines a homological intersection form $B: \H \times \H \to \Z$.  In the propositions and theorems that follow, we develop formulas that will help determine the parities of the intersection numbers for the homology classes discussed above.

\begin{notation}
The set of arrows $arr_{i,j}(\alpha)$ consists of the arrows with tail on $S_i$ and head on $S_j$.  For $e = (a,b) \in arr_{i,i}(\alpha)$, let $n(e) \in \Z$ denote
$$\#\{f \in arr_{i,i}(\alpha)\ |\ f \ \textnormal{links}\ e \ \textnormal{positively}\}
-\#\{f \in arr_{i,i}(\alpha)\ |\ f \ \textnormal{links}\ e \ \textnormal{negatively}\}.$$
Let $n_{i,j}(e)$ denote
$$\#\{f\in arr_{j,i}(\alpha)\ |\ \textnormal{head of}\ f\ \textnormal{in}\ (ab)^\circ \}
-\#\{f \in arr_{i,j}(\alpha)\ |\ \textnormal{tail of}\ f\ \textnormal{in}\ (ab)^\circ \},$$
and let $n^*_{i,j}(e)$ denote
$$\#\{f \in arr_{j,i}(\alpha)\ |\ \textnormal{head of}\ f\ \textnormal{in}\ (ba)^\circ \}
-\#\{f \in arr_{i,j}(\alpha)\ |\ \textnormal{tail of}\ f\ \textnormal{in}\ (ba)^\circ \}.$$
Notice that $n_{i,i}(e) = n(e)$ and $n^*_{i,i}(e) = -n(e)$.
\end{notation}

\begin{Prop}\label{P:BeS}
If $e = (a,b) \in arr_{i,i}(\alpha)$ and $i \neq j$, then
$$\begin{aligned}
B([e],[S_i]) &= B([e],[e]^*) = n(e) \\
B([e],[S_j]) &= n_{i,j}(e) \\
B([e]^*,[S_j]) &= n^*_{i,j}(e)
\end{aligned}$$
and
$$\begin{aligned}
B([e],[S_j]) &\equiv \#(p_i \cap o(v_j)) \mtwo\\
B([e]^*,[S_j]) &\equiv \#(p'_i \cap o(v_j)) \mtwo.
\end{aligned}$$
\end{Prop}
\begin{proof}
The proof of the first claim appears in \cite{Tvs}, but we include it here for completeness.  The loops $\rho(ab)$ and $\rho(ba)$ intersect transversely, except at their common origin $\rho(a)=\rho(b)$.  However, a small deformation makes these loops disjoint in a neighborhood of $\rho(a)=\rho(b)$.  Notice that there is a bijective correspondence between the transversal intersections and the arrows of $arr_{i,i}(\alpha)$ linked with $e$.  The sign of a transversal intersection with respect to $\rho(ab)$ is +1 when its corresponding arrow links $e$ positively, and -1 when its corresponding arrow links $e$ negatively.  Hence $B([e],[e]^*)=n(e)$.  Since $[S_i] = [e] + [e]^*$ in $\H$, it follows that
$$B([e],[S_i]) = B([e],[S_i]) - B([e],[e]) = B([e],[e]^*) = n(e).$$

Any intersections of the loops $\rho(ab)$ and $\rho(S'_j)$ are transversal intersections.  There is a bijective correspondence between these intersections and the union of the sets
$\{f\in arr_{j,i}(\alpha)\ |\ \textnormal{head of}\ f\ \textnormal{in}\ (ab)^\circ \}$ and 
$\{f \in arr_{i,j}(\alpha)\ |\ \textnormal{tail of}\ f\ \textnormal{in}\ (ab)^\circ \}$.
The sign of an intersection with respect to $\rho(ab)$ is +1 when its corresponding arrow is in the first set, and -1 when its corresponding arrow is in the second set.  Hence $B([e],[S_j]) = n_{i,j}(e)$.  Moreover, the arrows in the union of these two sets correspond to the letters that the sets $p_i$ and $o(v_j)$ have in common.  Hence $B([e],[S_j]) \equiv \#(p_i \cap o(v_j)) \mtwo$.  The claims about $[e]^*$ follow from similar arguments.
\end{proof}

\begin{Prop}\label{P:pSetsTonij}
If $i \in E$ is a double letter of the word $v_j$ in $p$, then
$\#(p_i) \equiv \sum_{k=1}^N n_{j,k}(e_i) \mtwo$ and
$\#(p'_i) \equiv \sum_{k=1}^N n^*_{j,k}(e_i) \mtwo$.
\end{Prop}
\begin{proof}
We have assumed that the letters in $p_i$ correspond to the arrows in $\alpha$ with exactly one endpoint on $(a_ib_i)^\circ$.  The other endpoint of such an arrow is either on $(b_ia_i)^\circ$ or on another component of $\alpha$.  For $k \neq j$,
$$\begin{aligned}
n_{j,k}(e_i) = \#\{f &\in arr_{k,j}(\alpha)\ |\ \textnormal{head of}\ f\ \textnormal{in}\ (a_ib_i)^\circ \}\\
	&-\#\{f \in arr_{j,k}(\alpha)\ |\ \textnormal{tail of}\ f\ \textnormal{in}\ (a_ib_i)^\circ \}\\
\equiv \#\{f &\in arr_{k,j}(\alpha)\ |\ \textnormal{head of}\ f\ \textnormal{in}\ (a_ib_i)^\circ \}\\
	&+\#\{f \in arr_{j,k}(\alpha)\ |\ \textnormal{tail of}\ f\ \textnormal{in}\ (a_ib_i)^\circ \} \mtwo,
\end{aligned}$$
but this sum is equal to the number of arrows between $S_j$ and $S_k$ with one endpoint on $(a_ib_i)^\circ$.  Notice
$$\begin{aligned}
n_{j,j}(e_i) = n(e_i) = &\#\{f \in arr_{j,j}(\alpha)\ |\ f \ \textnormal{links}\ e \ \textnormal{positively}\}\\
	&-\#\{f \in arr_{j,j}(\alpha)\ |\ f \ \textnormal{links}\ e \ \textnormal{negatively}\}\\
\equiv &\#\{f \in arr_{j,j}(\alpha)\ |\ f \ \textnormal{links}\ e \ \textnormal{positively}\}\\
	&+\#\{f \in arr_{j,j}(\alpha)\ |\ f \ \textnormal{links}\ e \ \textnormal{negatively}\} \mtwo,
\end{aligned}$$
with this sum equal to the number of arrows in $arr_{j,j}(\alpha)$ with one endpoint on $(a_ib_i)^\circ$ and the other on $(b_ia_i)^\circ$.  Hence $\#(p_i) \equiv \sum_{k=1}^{N} n_{j,k}(e_i)$.  The claim about $\#(p'_i)$ can be proved in a similar manner.
\end{proof}

\begin{Prop}\label{P:BSS}
For distinct core circles $S_i$ and $S_j$ of $\alpha$,
$$B([S_i],[S_j]) = \#(arr_{j,i}(\alpha)) - \#(arr_{i,j}(\alpha)).$$
The number of letters that the words $v_i$ and $v_j$ of $p$ have in common is equal to $\#(arr_{j,i}(\alpha)) + \#(arr_{i,j}(\alpha))$.  Therefore $B([S_i],[S_j])$ is equivalent modulo two to the number of letters that the words $v_i$ and $v_j$ of $p$ have in common.
\end{Prop}
\begin{proof}
Any intersections of the loops $\rho(S'_i)$ and $\rho(S'_j)$ are transverse intersections.  There is a bijective correspondence between these transverse intersections and the set of arrows between $S_i$ and $S_j$ in $\alpha$, that is the disjoint union $arr_{j,i}(\alpha) \cup arr_{i,j}(\alpha)$.  The sign of an intersection with respect to $\rho(S'_i)$ is +1 when its corresponding arrow is in $arr_{j,i}(\alpha)$, and -1 when its corresponding arrow is in $arr_{i,j}(\alpha)$.  Hence $B([S_i],[S_j]) = \#(arr_{j,i}(\alpha)) - \#(arr_{i,j}(\alpha))$.  The arrows in the disjoint union $arr_{j,i}(\alpha) \cup arr_{i,j}(\alpha)$ correspond to the letters that $v_i$ and $v_j$ have in common, so the total number of such letters is $\#(arr_{j,i}(\alpha)) + \#(arr_{i,j}(\alpha))$.
\end{proof}

\begin{Prop}\label{P:BeeNotInt}
If $i,j \in E$ are double letters of a word $w$ of $p$ that are not $w$-interlaced and $\#(w_i) \equiv \#(w_j) \equiv 0 \mtwo$, then $\#(w_i \cap w_j) \equiv B([e_i],[e_j]) \mtwo$.
\end{Prop}
\begin{proof}
The proof of this result is included as part the proof of Theorem 5.3.1 in \cite{Tvs}.  The same arguments apply because the intersections that contribute to $B([e_i],[e_j])$ are self-intersections of a single component of $\rho$.
\end{proof}

\begin{Prop}\label{P:BeeDiffWords}
Suppose $i \in E$ is a double letter of the word $v_m$ of $p$ and $j \in E$ is a double letter of word the $v_n$ of $p$, where $m \neq n$.  Then
$$\begin{aligned}
\#(p_i \cap p_j) &\equiv B([e_i],[e_j]) \mtwo \\
\#(p_i \cap p'_j) &\equiv B([e_i],[e_j]^*) \mtwo \\
\#(p'_i \cap p_j) &\equiv B([e_i]^*,[e_j]) \mtwo \\
\#(p'_i \cap p'_j) &\equiv B([e_i]^*,[e_j]^*) \mtwo
\end{aligned}$$
\end{Prop}
\begin{proof}
Any intersections of the loops $\rho(a_ib_i)$ and $\rho(a_jb_j)$ are transverse intersections.  There is a bijective correspondence between these transverse intersections and the set of arrows between $S_m$ and $S_n$ in $\alpha$ with one endpoint on $(a_ib_i)^\circ$ and the other on $(a_jb_j)^\circ$.  Such arrows correspond to the letters of $p$ appearing in $p_i \cap p_j$.  Hence $\#(p_i \cap p_j) \equiv B([e_i],[e_j]) \mtwo$.  The other three claims follow from similar arguments.
\end{proof}

\begin{Thm}\label{T:SnDpHomology}
For all $1 \leq n \leq N$ and $d \in \D_p$,
$$B([S_n],[C_d]) \equiv W(v_n,d) \mtwo.$$
\end{Thm}
\begin{proof}
If $d(n)$ is an empty sequence, then the intersections of $\rho(S'_n)$ and $\rho(C'_d)$ are all transversal intersections.  These intersections are in bijective correspondence with the arrows of $\alpha$ with one endpoint on $S_n$ and the other on any arc $C_d(k)$, with $k \neq n$ since $d(n) = \emptyset$.  If an arrow of $\alpha$ has one endpoint on $S_n$ and the other on $C_d(k)$, then it corresponds to a letter of $p$ that $v_n$ and $d(k)$ have in common.  Hence $B([S_n],[C_d])$ is equivalent to $\sum_{k \neq n} \#(o(v_n) \cap o(d(k))) = W(v_n,d)$ modulo two in this case.

Now assume $d(n)$ is a nonempty sequence with first letter $i$ and last letter $j$.  Then $d(n)$ corresponds to an arc $C'_d(n)$ on $S'_n$.  A neighborhood of the segment $\rho(C'_d(n))$ in $\Sigma$ is depicted in Figure \ref{F:SnCase}, where the thickened horizontal line segment represents the segment of $\rho(S'_n)$ between crossings $i$ and $j$.  The vertical lines depict neighborhoods of crossings $i$ and $j$ in some components $\rho(S'_t)$ and $\rho(S'_u)$, respectively, and their ends are labeled with the letters $i_1,i_2$ and $j_1,j_2$ to represent the letters flanking $i$ and $j$ in the words $v_t$ and $v_u$ of $p$, respectively.

An intersection of the loops $\rho(S'_n)$ and $\rho(C'_d)$ is either a transversal intersection involving $\rho(S'_n)$ and a segment of $\rho(C'_d)$ not on $\rho(S'_n)$ or an intersection resulting from the segment $\rho(C'_d(n))$ on $\rho(S'_n)$.  As in the above discussion, the transversal intersections correspond to letters of $p$ that $v_n$ and any $d(k)$ for $k \neq n$ have in common, so there are a total of $\sum_{k \neq n} \#(o(v_n) \cap o(d(k)))$ such intersections.  The other intersections require further consideration.

If the points $i,i_1,i_2$ lie in the cyclic order $i_1ii_2$ (respectively $i_2ii_1$) on $S'_t$, then the sign of the crossing $i$ with respect to $\rho(S'_n)$ is -1 (respectively +1) and $i \in A_n \cup A'_n$ (respectively $i \in A_t \cup A'_t$).  Similarly, if the points $j,j_1,j_2$ lie in the cyclic order $j_1jj_2$ (respectively $j_2jj_1$) on $S'_u$, then the sign of the crossing $j$ with respect to $\rho(S'_n)$ is -1 (respectively +1) and $j \in A_n \cup A'_n$ (respectively $j \in A_u \cup A'_u$).  The possible cyclic orders of the points $i,i_1,i_2$ and $j,j_1,j_2$ on $S'_t$ and $S'_u$, respectively, yield the following four cases:
\begin{enumerate}
\item{$i_1ii_2$ on $S'_t$ and $j_1jj_2$ on $S'_u$;}
\item{$i_1ii_2$ on $S'_t$ and $j_2jj_1$ on $S'_u$;}
\item{$i_2ii_1$ on $S'_t$ and $j_1jj_2$ on $S'_u$;}
\item{$i_2ii_1$ on $S'_t$ and $j_2jj_1$ on $S'_u$.}
\end{enumerate}

Consider the loop $\sigma$ in $\Sigma$ that is the ``push off" of $\rho(S'_n)$ to the right in a neighborhood of $\rho(S'_n)$.  In particular, for the above four cases, consider a neighborhood of the segment $\rho(C'_d(n))$, as depicted in Figure \ref{F:SnCase}.  The dashed lines in the figure represent $\sigma$ and the thickened lines represent a portion of the loop $\rho(C'_d)$.  The intersections resulting from segment $\rho(C'_d(n))$ on $\rho(S'_n)$ are transversal intersections of $\sigma$ and $\rho(C'_d)$ in this neighborhood.  In every case, the number of intersections on the interior of the segment $\rho(C'_d(n))$ is the same.  In a neighborhood of a self-intersection of $\rho(S'_n) \cap \rho(C'_d(n))$, there are two crossings involving $\sigma$ and the interior of $\rho(C'_d(n))$.  These two crossings have opposite signs, so they do not contribute to $B([S_n],[C_d])$.  Therefore we only consider those intersections of $\sigma$ and $\rho(C'_d(n))$ that do not result from self-intersections of $\rho(S'_n) \cap \rho(C'_d(n))$.  Such intersections bijectively correspond to letters that the subsequences $x_1$ and $x_2$ of $v_n = ix_1jx_2$ have in common, and there are $\#(o(ix_1j) \cap o(jx_2i))$ such letters.  Notice that there may be additional transverse intersections of $\sigma$ and $\rho(C'_d)$ in each of the four cases, resulting from the portions of the remaining segments of $\rho(C'_d)$ depicted in Figure \ref{F:SnCase}.  Specifically, there are 1, 0, 2, and 1 additional intersections in cases (1), (2), (3), and (4), respectively, and these totals are equivalent to $\gamma_{n}^P(i,j)$ modulo two in each case.  Hence $B([S_n],[C_d])$ is equivalent to
$$\#(o(ix_1j) \cap o(jx_2i)) + \gamma_{n}^{P}(i,j) + \sum_{k \neq n} \# (o(v_n) \cap o(d(k))) = W(v_n,d)$$
modulo two.
\end{proof}

\begin{figure}
\scalebox{.32}{\includegraphics{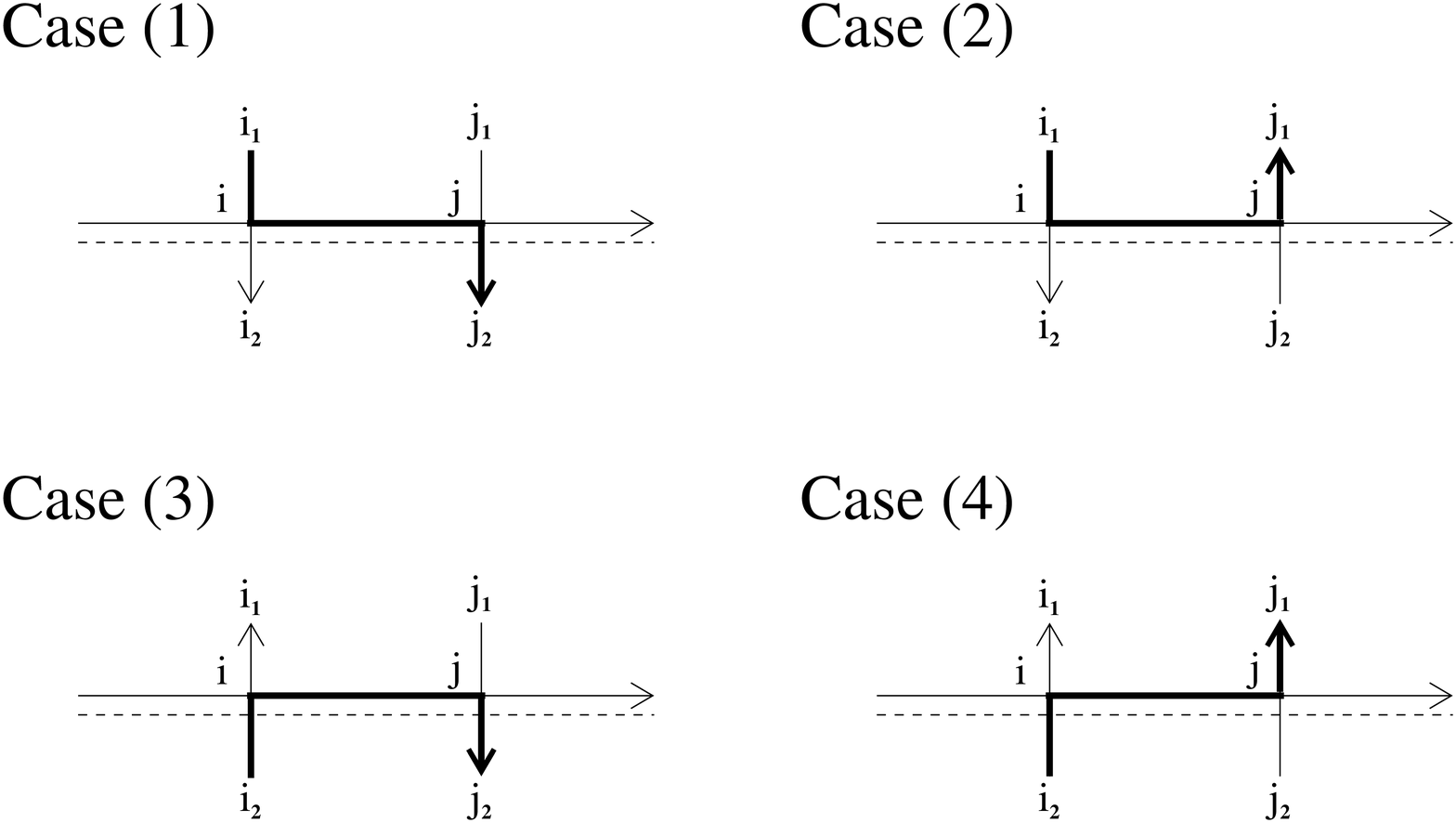}}
\caption{Neighborhood of the segment $\rho(C'_d(n))$.}\label{F:SnCase}
\end{figure}

\begin{Thm}\label{T:eiDpHomology}
If a word $v_n$ of $p$ has a double letter $i \in E$, then
$$B([e_i],[C_d]) \equiv Q_n(p_i,d) \mtwo$$
and
$$B([e_i]^*,[C_d]) \equiv Q_n(p'_i,d) \mtwo.$$
\end{Thm}
\begin{proof}
Since $i$ is a double letter of the word $v_n$, it follows that $v_n = ix_1ix_2$ for some subsequences $x_1,x_2$.  Recall that we have assumed $p_i$ and $p'_i$ contain the letters corresponding to the arrows in $\alpha$ with exactly one endpoint on $(a_ib_i)^\circ$ and $(b_ia_i)^\circ$, respectively.  The letters in $p_i$ then correspond to the letters that occur exactly once in either $x_1$ or $x_2$, but we may assume it is $x_1$.  Then the letters in $p'_i$ correspond to the letters that occur exactly once in $x_2$.  

If $d(n)$ is an empty sequence, then the intersections of the loops $\rho(a_ib_i)$ and $\rho(C'_d)$ are all transversal intersections.  These intersections are in bijective correspondence with the arrows of $\alpha$ with one endpoint on $(a_ib_i)^\circ$ and the other on any arc $C_d(k)$, with $k \neq n$ since $d(n) = \emptyset$.  If an arrow of $\alpha$ has one endpoint on $(a_ib_i)^\circ$ and the other on $C_d(k)$, then it corresponds to a letter of $p$ that $p_i$ and $o(d(k))$ have in common.  Therefore $B([e_i],[C_d]) \equiv \sum_{k \neq n}\#(p_i \cap o(d(k))) \mtwo$.  Hence $B([e_i],[C_d]) \equiv Q_n(p_i,d) \mtwo$ in this case since $\delta_n(ix_1i,d(n)) = \delta_n(ix_1i,\emptyset) = 0$ and $\epsilon_n^P(ix_1i,d(n)) = \epsilon_n^P(ix_1i,\emptyset) = 0$.

Now suppose $d(n) = y x z$ for some single letters $y,z$ of $v_n$ and subsequence $x$ of $v_n$.  The subsequence $d(n)$ of $v_n$ corresponds to an arc $C'_d(n) \subset S'_n$ with initial point $y$ and terminal point $z$.  The subsequence $ix_1i$ of $v_n$ corresponds to the arc $a_ib_i \subset S'_n$ with initial point $a_i$ and terminal point $b_i$.  Notice that the points $a_i,b_i,y,z$ must be pairwise distinct since they are the four endpoints of two different arrows in $\alpha$.

An intersection of $\rho(a_ib_i)$ and $\rho(C'_d)$ is either a transversal intersection involving $\rho(a_ib_i)$ and a segment of $\rho(C'_d)$ not on $\rho(S'_n)$ or an intersection resulting from the segment $\rho(C'_d(n))$ on $\rho(S'_n)$.  As in the above discussion, the transversal intersections correspond to letters of $p$ that $p_i$ and $o(d(k))$ for $k \neq n$ have in common, so there are a total of $\sum_{k \neq n} \#(p_i \cap o(d(k)))$ such intersections.

Now consider the intersections resulting from the segment $\rho(C'_d(n))$ on $\rho(S'_n)$.  The six possible cyclic orders of the points $a_i,b_i,y,z$ on $S'_n$ are:
\begin{enumerate}
\item{$a_i, y, b_i, z$;}
\item{$a_i, z, b_i, y$;}
\item{$a_i, y, z, b_i$;}
\item{$a_i, z, y, b_i$;}
\item{$a_i, b_i, y, z$;}
\item{$a_i, b_i, z, y$.}
\end{enumerate}
Consider the loop $\sigma$ in $\Sigma$ that is the ``push off" of $\rho(a_ib_i)$ to the right in a neighborhood of $\rho(S'_n)$.  The intersections resulting from the segment $\rho(C'_d(n))$ are transversal intersections of $\sigma$ and $\rho(C'_d)$ in this neighborhood.  In particular, for the above six cases, consider a neighborhood of $\rho(S'_n)$ in $\Sigma$ as depicted in Figure \ref{F:eiCase}.  The dashed loops in the figure represent $\sigma$ and the thickened lines represent the segment $\rho(C'_d(n))$.  The segments of $\rho(C'_d)$ immediately before and after $\rho(C'_d(n))$ are not indicated because they depend on the signs of crossings $y$ and $z$ with respect to $\rho(S'_n)$.  If the sign of crossing $y$ is -1 (respectively +1), then $y \in A_n \cup A'_n$ (respectively $y \notin A_n \cup A'_n$).  If the sign of crossing $z$ is -1 (respectively +1), then $z \in A_n \cup A'_n$ (respectively $z \notin A_n \cup A'_n$).

\begin{figure}
\scalebox{.25}{\includegraphics{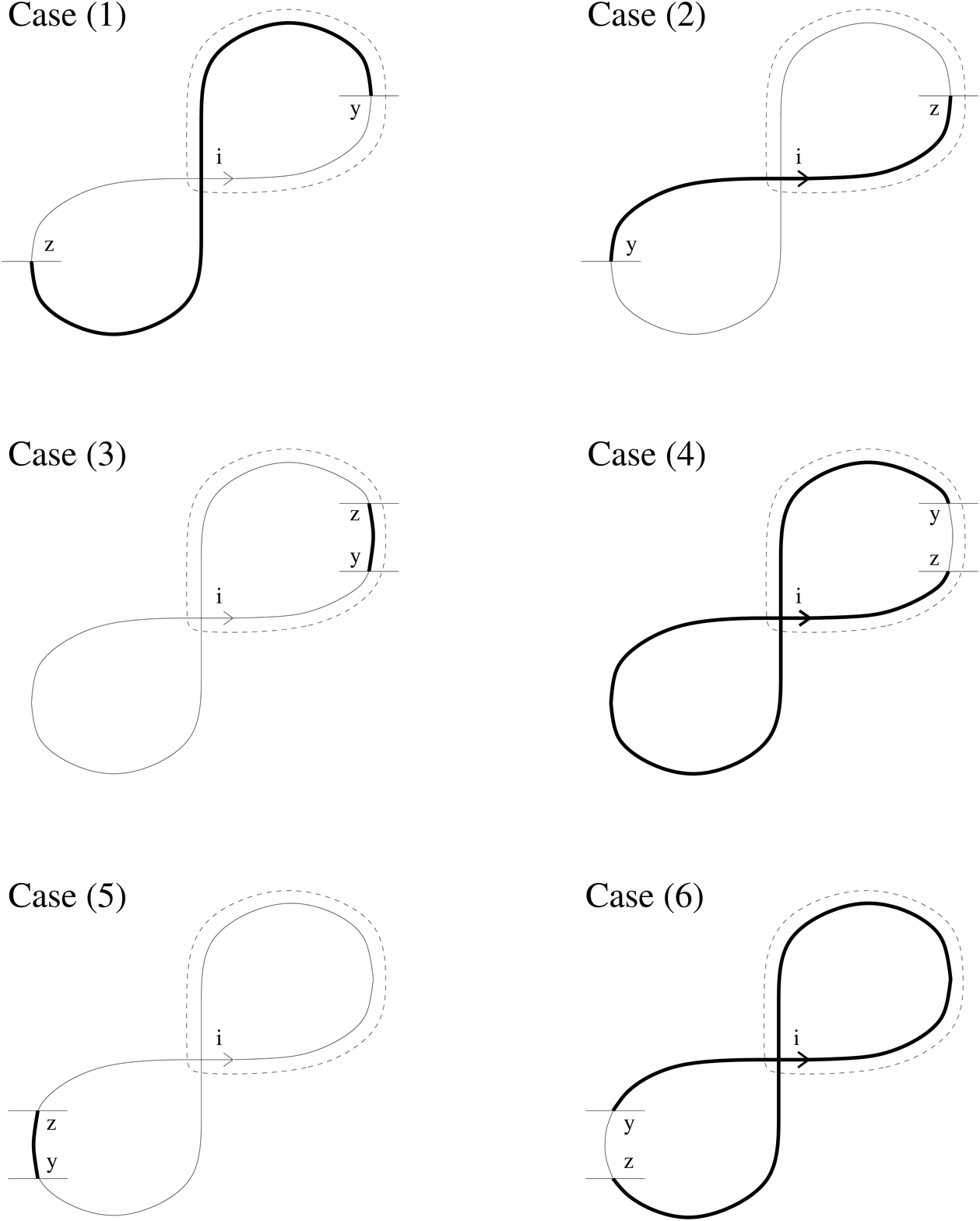}}
\caption{Neighborhood of the loop $\rho(S'_n)$.}\label{F:eiCase}
\end{figure}

In each of the six cases, there may be transversal intersections of $\sigma$ and $\rho(C'_d)$ in neighborhoods of crossings $i$, $y$, and $z$.  In a neighborhood of crossing $i$, there are two intersections of $\sigma$ and $\rho(C'_d(n))$ in Cases (4) and (6), and there is one intersection in Cases (1) and (2).  In a neighborhood of crossing $y$, there is an intersection of $\sigma$ and the segment of $\rho(C'_d)$ immediately before $\rho(C'_d(n))$ in Cases (1), (3), and (4) if and only if the sign of crossing $y$ with respect to $\rho(S'_n)$ is +1.  In a neighborhood of crossing $z$, there is an intersection of $\sigma$ and the segment of $\rho(C'_d)$ immediately after $\rho(C'_d(n))$ in Cases (2), (3), and (4) if and only if the sign of crossing $z$ with respect to $\rho(S'_n)$ is -1.  Observe that there are no intersections of $\sigma$ and $\rho(C'_d)$ in these three neighborhoods otherwise.  Since $p_i$ contains the letters that occur exactly once in the subsequence $x_1$ of $v_n$, it follows that the total number of transversal intersections of $\sigma$ and $\rho(C'_d)$ in these neighborhoods of crossings $i$, $y$, and $z$ is given by $\epsilon_n^P(ix_1i,d(n))$.

The only remaining transversal intersections of $\sigma$ and $\rho(C'_d)$ left to consider are intersections of $\sigma$ and $\rho(C'_d(n))$ that do not occur in the three crossing neighborhoods discussed above.  In a neighborhood of a self-intersection of $\rho(a_ib_i) \cap \rho(C'_d(n))$, there are two crossings involving $\sigma$ and the interior of $\rho(C'_d(n))$.  These two crossings have opposite signs, so they do not contribute to $B([e_i],[C_d])$.  Therefore we only consider those intersections of $\sigma$ and $\rho(C'_d(n))$ that do not result from self-intersections of $\rho(a_ib_i) \cap \rho(C'_d(n))$.  Such intersections bijectively correspond to arrows in $arr_{n,n}(\alpha)$ with one endpoint on arc $a_ib_i$ and the other on $C_d(n)$, provided both endpoints are not on $a_ib_i \cap C_d(n)$.  In the following arguments, let the labels $s_1,s_2,s_3,s_4$ denote subsequences of the word $v_n$ and assume the endpoints $a_i$ and $b_i$ correspond to the first and second instances of $i$, respectively.  In Case (1), if $v_n = i\ s_1\ y\ s_2\ i\ s_3\ z\ s_4$, then these arrows correspond to letters that $s_1$ and $s_2$ have in common, $s_1$ and $s_3$ have in common, or $s_2$ and $s_3$ have in common.   In Case (2), if $v_n = i\ s_1\ z\ s_2\ i\ s_3\ y\ s_4$, then these arrows correspond to letters that $s_1$ and $s_4$ have in common, $s_2$ and $s_4$ have in common, or $s_2$ and $s_1$ have in common.  In Case (3), if $v_n = i\ s_1\ y\ s_2\ z\ s_3\ i\ s_4$, then these arrows correspond to letters that subsequences $s_1$ and $s_2$ have in common or letters that $s_3$ and $s_2$ have in common.  In Case (4), if $v_n = i\ s_1\ z\ s_2\ y\ s_3\ i\ s_4$, then these arrows correspond to letters that $s_1$ and $s_4$ have in common, $s_2$ and $s_3$ have in common, $s_2$ and $s_4$ have in common, $s_2$ and $s_1$ have in common, or $s_3$ and $s_4$ have in common.  In Case (5), if $v_n = i\ s_1\ i\ s_2\ y\ s_3\ z\ s_4$, then these arrows correspond to letters that subsequences $s_1$ and $s_3$ have in common.  In Case (6), if $v_n = i\ s_1\ i\ s_2\ z\ s_3\ y\ s_4$, then these arrows correspond to letters that $s_1$ and $s_4$ have in common or $s_1$ and $s_2$ have in common.  Hence the total number of such arrows is $\delta_n(ix_1i,d(n))$ in each of the six cases.

Having considered all possible intersections that contribute to $B([e_i],[C_d])$, we conclude that $B([e_i],[C_d])$ is equivalent to
$$\delta_n(ix_1i,d(n)) + \epsilon_n^P(ix_1i,d(n)) + \sum_{k \neq n}\#(p_i \cap o(d(k))) = Q_n(p_i,d)$$
modulo two.  The claim $B([e_i]^*,[C_d]) \equiv Q_n(p'_i,d) \mtwo$ can be proven with similar arguments.
\end{proof}

\begin{Thm}\label{T:DpDpHomology}
If $d_1,d_2 \in \D$, then 
$$B([C_{d_1}],[C_{d_2}]) \equiv \sum_{n=1}^N D_n(d_1,d_2) \mtwo.$$
\end{Thm}
\begin{proof}
If at least one of $d_1(n)$ and $d_2(n)$ is empty for every $1 \leq n \leq N$, then all intersections of $\rho(C'_{d_1})$ and $\rho(C'_{d_2})$ are transversal intersections since these loops have segments on distinct components of $\rho$.  Suppose $d_1(n) = a c b$ for some single letters $a,b$ of $v_n$ and subsequence $c$ of $v_n$.  The subsequence $d_1(n)$ of $v_n$ corresponds to an arc $C'_{d_1}(n)$ in $S'_n$ with initial point $a$ and terminal point $b$.  The transversal intersections of the arc $\rho(C'_{d_1}(n))$ and $\rho(C'_{d_2})$ are in bijective correspondence with the arrows of $\alpha$ with one endpoint on arc $C_{d_1}(n)$ and the other on any arc $C_{d_2}(k)$, with $k \neq n$ since $d_2(n) = \emptyset$.  If an arrow of $\alpha$ has one endpoint on $C_{d_1}(n)$ and the other on $C_{d_2}(k)$, then it corresponds to a letter of $p$ that $d_1(n)$ and $d_2(k)$ have in common.  Therefore the total number of intersections of $\rho(C'_{d_1}(n))$ and $\rho(C'_{d_2})$ is $\sum_{k \neq n} \#(o(d_1(n)) \cap o(d_2(k)))$, with $\#(o(d_1(n)) \cap o(d_2(k))) = 0$ for any $k$ with $d_2(k) = \emptyset$.  Hence $B([C_{d_1}],[C_{d_2}]) \equiv \sum_{n=1}^N D_n(d_1,d_2) \mtwo$ since $\delta_n(d_1(n),d_2(n)) = \epsilon_n^P(d_1(n),d_2(n)) = \sum_{k \neq n} \#(o(d_1(n)) \cap o(d_2(k))) = 0$ when $d_1(n)$ is an empty sequence and $\delta_n(d_1(n),d_2(n)) = \epsilon_n^P(d_1(n),d_2(n)) = 0$ when $d_2(n)$ is an empty sequence.

Now remove the hypothesis that at least one of $d_1(n)$ and $d_2(n)$ is empty for every $1 \leq n \leq N$.  Again suppose $d_1(n) = a c b$ for some single letters $a,b$ of $v_n$ and subsequence $c$ of $v_n$.  If $d_2(n) = \emptyset$, then the total number of intersections of $\rho(C'_{d_1}(n))$ and $\rho(C'_{d_2})$ is $\sum_{k \neq n} \#(o(d_1(n)) \cap o(d_2(k))) = D_n(d_1,d_2)$ as shown above since $\delta_n(d_1(n),\emptyset) = \epsilon_n^P(d_1(n),\emptyset) = 0$.
  
Now suppose $d_2(n) = y x z$ for some single letters $y,z$ of $v_n$ and subsequence $x$ of $v_n$.  For any $k \neq n$, the intersections of the segments $\rho(C'_{d_1}(n))$ and $\rho(C'_{d_2}(k))$ are transversal intersections.  These intersections are in bijective correspondence with the arrows of $\alpha$ with one endpoint on arc $C_{d_1}(n)$ and the other on arc $C_{d_2}(k)$.  Such an arrow corresponds to a letter of $p$ that $d_1(n)$ and $d_2(k)$ have in common.  Therefore the total number of intersections of $\rho(C'_{d_1}(n))$ and the segments of $\rho(C'_{d_2})$ not on $\rho(S'_n)$ is $\sum_{k \neq n} \#(o(d_1(n)) \cap o(d_2(k)))$.

The intersections of $\rho(C'_{d_1})$ and $\rho(C'_{d_2})$ on $\rho(S'_n)$ may also contribute to \linebreak $B([C_{d_1}],[C_{d_2}])$.  Consider the loop $\sigma$ in $\Sigma$ that is the ``push off" of $\rho(C'_{d_1})$ to the right in a neighborhood of $\rho$.  The intersections of $\rho(C'_{d_1})$ and $\rho(C'_{d_2})$ on $\rho(S'_n)$ are transversal intersections of $\sigma$ and $\rho(C'_{d_2})$ in a neighborhood of $\rho(S'_n)$.

The six possible cyclic orders of the points $a,b,y,z$ on $S'_n$ are:
\begin{enumerate}
\item{$a, y, b, z$;}
\item{$a, z, b, y$;}
\item{$a, y, z, b$;}
\item{$a, z, y, b$;}
\item{$a, b, y, z$;}
\item{$a, b, z, y$,}
\end{enumerate}
with the possibility that $a,y,z$ and $b,y,z$ are not pairwise distinct where appropriate (but $a,b$ and $y,z$ are distinct).  For these six cases, consider a neighborhood of the loop $\rho(S'_n)$ in $\Sigma$ as depicted in Figure \ref{F:DpCase}.  The dashed segments in the figure represent segments of $\sigma$ and the thickened lines represent the segment $\rho(C'_{d_2}(n))$.  The segments of $\sigma$ immediately before and after $\rho(C'_{d_1}(n))$ are not indicated because they depend on the signs of crossings $a$ and $b$ with respect to $\rho(S'_n)$.  The segments of $\rho(C'_{d_2})$ immediately before and after $\rho(C'_{d_2}(n))$ are also not indicated because they depend on the signs of crossings $y$ and $z$ with respect to $\rho(S'_n)$.  If the sign of one of the crossings is -1 (respectively +1), then the letter it is labeled by appears in $A_n \cup A'_n$ (respectively does not appear in $A_n \cup A'_n$).

\begin{figure}
\scalebox{.25}{\includegraphics{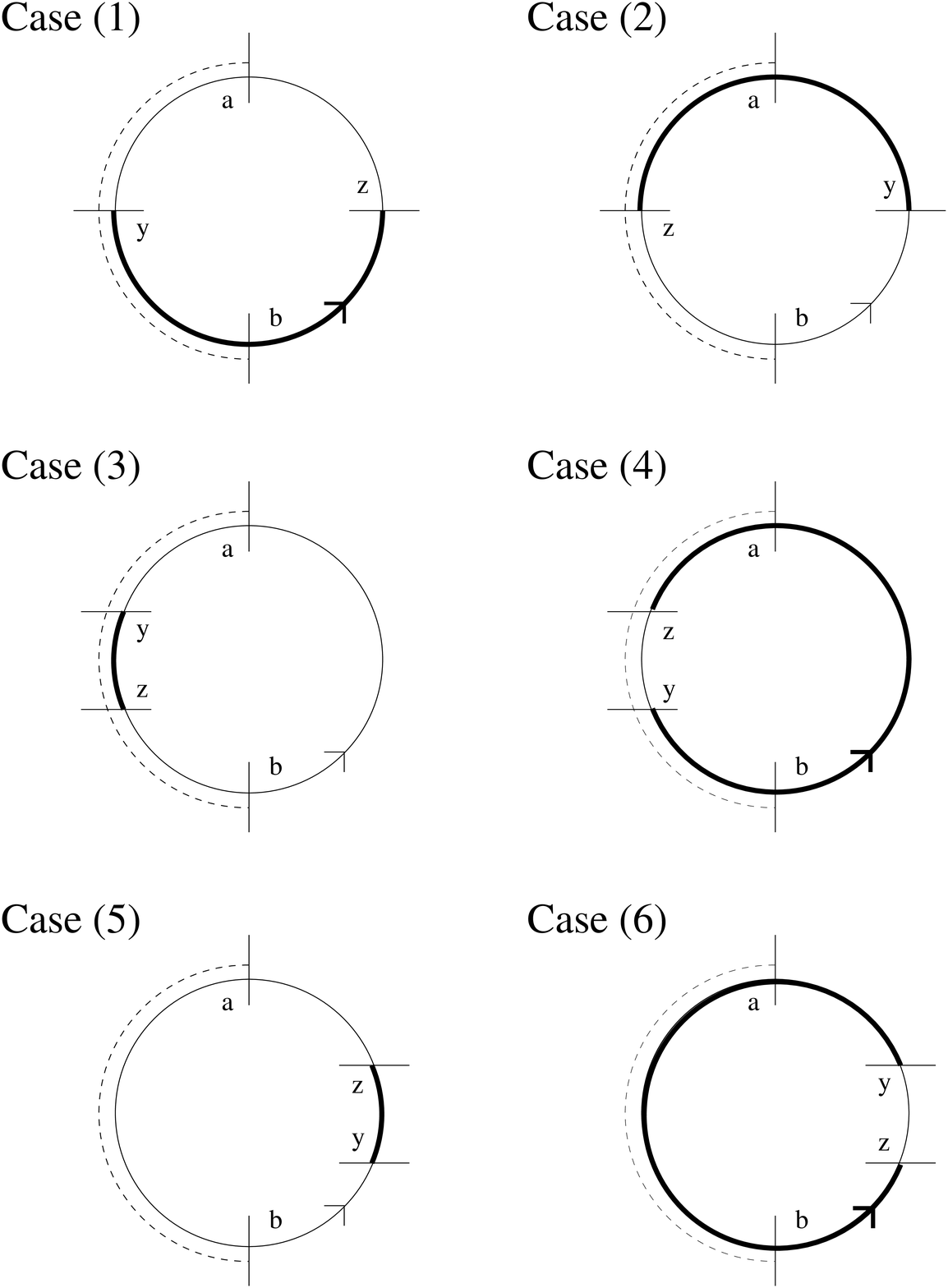}}
\caption{Neighborhood of the loop $\rho(S'_n)$.}\label{F:DpCase}
\end{figure}

In each of the six cases, there may be transversal intersections of $\sigma$ and $\rho(C'_{d_2})$ in neighborhoods of crossings $a$, $b$, $y$, and $z$.  If the sign of crossing $a$ with respect to $\rho(S'_n)$ is -1 and $a,y,z$ are pairwise distinct, then there is one intersection of $\sigma$ and $\rho(C'_{d_2}(n))$ in a neighborhood of crossing $a$ in cases (2), (4), and (6).  It is possible for $a=y$ in Cases (1), (2), (3), and (6), and if so there are no intersections of $\sigma$ and $\rho(C'_{d_2})$ in a neighborhood of crossing $a=y$.  It is possible for $a=z$ in Cases (1), (2), (4), and (5), and if so there are two intersections of $\sigma$ and $\rho(C'_{d_2})$ in a neighborhood of crossing $a=z$ when the sign of crossing $a=z$ is -1 and none when the sign of crossing $a=z$ is +1.  Observe that there are no intersections of $\sigma$ and $\rho(C'_{d_2})$ in a neighborhood of crossing $a$ otherwise.

If the sign of crossing $b$ with respect to $\rho(S'_n)$ is +1 and $b,y,z$ are pairwise distinct, then there is one intersection of $\sigma$ and $\rho(C'_{d_2}(n))$ in a neighborhood of crossing $b$ in cases (1), (4), and (6).  It is possible for $b=y$ in Cases (1), (2), (4), and (5), and if so there are two intersections of $\sigma$ and $\rho(C'_{d_2})$ in a neighborhood of crossing $b=y$ when the sign of crossing $b=y$ is +1 and none when the sign of crossing $b=y$ is -1.  It is possible for $b=z$ in Cases (1), (2), (3), and (6), and if so there are no intersections of $\sigma$ and $\rho(C'_{d_2})$ in a neighborhood of crossing $b=z$.  Observe that there are no intersections of $\sigma$ and $\rho(C'_{d_2})$ in a neighborhood of crossing $b$ otherwise.

Assume $a,b,y,z$ are pairwise distinct.  If the sign of crossing $y$ with respect to $\rho(S'_n)$ is +1, then there is one intersection of $\sigma$ and $\rho(C'_{d_2})$ in a neighborhood of crossing $y$ in cases (1), (3), and (4).  If the sign of crossing $z$ with respect to $\rho(S'_n)$ is -1, then there is one intersection of $\sigma$ and $\rho(C'_{d_2})$ in a neighborhood of crossing $z$ in cases (2), (3), and (4).  Observe that there are no intersections of $\sigma$ and $\rho(C'_{d_2})$ in these neighborhoods of crossings $y$ and $z$ otherwise when $a,b,y,z$ are pairwise distinct.

Combining the intersection information in the previous three paragraphs, notice that the total number of transversal intersections of $\sigma$ and $\rho(C'_{d_2})$ in neighborhoods of crossings $a$, $b$, $y$, and $z$ is given by $\epsilon_n^P(d_1(n),d_2(n))$.

The only remaining transversal intersections of $\sigma$ and $\rho(C'_{d_2})$ in a neighborhood of $\rho(S'_n)$ in $\Sigma$ to consider are intersections of $\sigma$ and $\rho(C'_d(n))$ that do not occur in the four crossing neighborhoods discussed above.  In a neighborhood of a self-intersection of $\rho(C'_{d_1}(n)) \cap \rho(C'_{d_2}(n))$, there are two crossings involving $\sigma$ and the interior of $\rho(C'_{d_2}(n))$.  These two crossings have opposite signs, so they do not contribute to $B([C_{d_1}],[C_{d_2}])$.  Therefore we only consider those intersections of $\sigma$ and $\rho(C'_{d_2}(n))$ that do not result from self-intersections of $\rho(C'_{d_1}(n)) \cap \rho(C'_{d_2}(n))$.  Such intersections bijectively correspond to arrows in $arr_{n,n}(\alpha)$ with one endpoint on arc $C_{d_1}(n)$ and the other on $C_{d_2}(n)$, provided both endpoints are not on the arc $C_{d_1}(n) \cap C_{d_2}(n)$.  In the following arguments, let the labels $s_1,s_2,s_3,s_4$ denote subsequences of the word $v_n$.  In Case (1), if $v_n = a\ s_1\ y\ s_2\ b\ s_3\ z\ s_4$, then these arrows correspond to letters that $s_1$ and $s_2$ have in common, $s_1$ and $s_3$ have in common, or $s_2$ and $s_3$ have in common.   In Case (2), if $v_n = a\ s_1\ z\ s_2\ b\ s_3\ y\ s_4$, then these arrows correspond to letters that $s_1$ and $s_4$ have in common, $s_2$ and $s_4$ have in common, or $s_2$ and $s_1$ have in common.  In Case (3), if $v_n = a\ s_1\ y\ s_2\ z\ s_3\ b\ s_4$, then these arrows correspond to letters that subsequences $s_1$ and $s_2$ have in common or letters that $s_3$ and $s_2$ have in common.  In Case (4), if $v_n = a\ s_1\ z\ s_2\ y\ s_3\ b\ s_4$, then these arrows correspond to letters that $s_1$ and $s_4$ have in common, $s_2$ and $s_3$ have in common, $s_2$ and $s_4$ have in common, $s_2$ and $s_1$ have in common, or $s_3$ and $s_4$ have in common.  In Case (5), if $v_n = a\ s_1\ b\ s_2\ y\ s_3\ z\ s_4$, then these arrows correspond to letters that subsequences $s_1$ and $s_3$ have in common.  In Case (6), if $v_n = a\ s_1\ b\ s_2\ z\ s_3\ y\ s_4$, then these arrows correspond to letters that $s_1$ and $s_4$ have in common or $s_1$ and $s_2$ have in common.  Hence the total number of such arrows is $\delta_n(d_1(n),d_2(n))$ in each of the six cases.

Having considered all possible intersections of segment $\rho(C'_{d_1}(n))$ and $\rho(C'_{d_2})$, we conclude that there are a total of
$$\delta_n(d_1(n),d_2(n)) + \epsilon_n^P(d_1(n),d_2(n)) + \sum_{k \neq n}\#(o(d_1(n)) \cap o(d_2(k)))$$
such intersections.  Thus $B([C_{d_1}],[C_{d_2}]) \equiv \sum_{n=1}^N D_n(d_1,d_2) \mtwo$.
\end{proof}
\end{section}

\begin{section}{Statement of the Main Theorem}
\begin{MainThm}
A pair (a Gauss paragraph $p$ in an alphabet $E$, a word-wise partition $P$ of $E$) is realizable by a closed curve on $S^2$ if and only if the following seven conditions are satisfied:
\begin{enumeratei}
    \item{if $i \in E$ is a double letter of a word $w$ of $p$, then
         $$\#(w_i) \equiv \#(p_i) \equiv \#(p'_{i}) \equiv 0 \mtwo;$$}
    \item{if $i \in E$ is a double letter of a word of $p$ that is not the word $w$, then
         $$\#(p_i \cap o(w)) \equiv \#(p'_i \cap o(w)) \equiv 0 \mtwo;$$}
    \item{any two words in $p$ have an even number of letters in common;}
    \item{if $i,j \in E$ are both double letters of the word $w$ of $p$ and are not
         $w$-interlaced, then $\#(w_i \cap w_j) \equiv 0 \mtwo$;}
    \item{$P$ is compatible with $p$;}
    \item{if $i \in E$ is a double letter of a word of $p$ and $j \in E$ is a double letter of a different word of $p$, then
         $$\#(p_i \cap p_j) \equiv \#(p_i \cap p'_j) \equiv \#(p'_i \cap p_j) \equiv \#(p'_i \cap p'_j) \equiv 0 \mtwo;\ and$$}
    \item{$P$ is compatible with $\D_p$.}
\end{enumeratei}
\end{MainThm}
Two corollaries of this theorem are Theorem 5.3.1 of Turaev's paper \cite{Tvs} and the central theorem of Rosenstiehl's paper \cite{Rsad}.  Condition (i) is analogous to Gauss' famous necessary condition for closed curves with a single component in \cite{Gw}, and Condition (iv) is analogous to the second condition of Rosenstiehl's central theorem in \cite{Rsad}.  The theorem above will be proved in the next section.  The following lemmas will be used in the proof.

\begin{Lemma}\label{L:LengthsEven}
Suppose any two words of a Gauss paragraph $p$ have an even number of letters in common.  Then every word of $p$ has even length.
\end{Lemma}
\begin{proof}
Recall that the length of a word $w$ of $p$ is
$$2(\#\{\textnormal{double letters in}\ w\}) + \#\{\textnormal{single letters in}\ w\}.$$
Therefore the length of $w$ is equivalent modulo two to the number of single letters of $w$.  But the number of single letters in $w$ is the sum of the numbers of letters that $w$ has in common with each other word of $p$.  Since each term of this sum is even, it follows that $w$ has even length. 
\end{proof}

\begin{Lemma}\label{L:CondVString}
If a pair (a Gauss paragraph $p$ in an alphabet $E$, a word-wise partition $P$ of $E$) satisfies conditions (i), (iii), and (v) of the main theorem, then a virtual string $\alpha$ can be constructed such that $\alpha$ has underlying Gauss word $p$ and induced word-wise partition $P$.
\end{Lemma}
\begin{proof}
The three conditions allow us to construct a virtual string $\alpha$ from $(p,P)$ in the following manner.  Suppose $p = (v_1,v_2,...,v_N)$ and $P = (A_1 \cup A'_1,A_2 \cup A'_2,...,A_N \cup A'_N)$.  Each word $v_n$ of $p$ can be written in the form $z_1z_2 \cdots z_{2M_n}$ for some $M_n \in \mathbb{N}$ since each word of $p$ has even length by Condition (iii) and Lemma \ref{L:LengthsEven}.  Consider a copy $R_n$ of $\R$ with distinguished points $\{1,2,...,2M_n\}$.  The core circle $S_n$ of $\alpha$ corresponding to $v_n$ will be the union of $R_n$ and a point at infinity with orientation extending the right-handed orientation on $R_n$.  A distinguished point $a$ of this core circle will correspond to the endpoint of an arrow in $\alpha$ for the letter $z_a \in E$.  If word $v_n$ has a double letter, then $z_a = z_b$ for some indices $a$ and $b$.  Note $a-b \equiv 1 \mtwo$ by Condition (i), so one of the indices $a,b$ is odd and the other is even.  A consequence of Condition (v) is that for single letters $z_x,z_y$ of $v_n$ in $A_n \cup A'_n$, we have $x-y \equiv 0 \mtwo$ if $z_x,z_y$ appear in the same set of this union and $x-y \equiv 1 \mtwo$ if $z_x,z_y$ appear in different sets of this union.  Therefore the single letters of $A_n$ are all in either even or odd positions in $z_1z_2 \cdots z_{2M_n}$, and the single letters of $A'_n$ all have the other position parity.

Any two distinct words $v_j, v_k$ of $p$ have an even number $2m$ (with $m \in \N$) of letters of common by Condition (iii), $m$ of which appear in $A_j \cup A'_j$ and $m$ of which appear in $A_k \cup A'_k$ by the definition of word-wise partition.  If a single letter of both $v_j$ and $v_k$ is an element of $A_j \cup A'_j$, then attach an arrow to the points in $R_j$ and $R_k$ corresponding to the letter, with the arrow's tail on $R_j$ and its head on $R_k$.  If a single letter of both $v_j$ and $v_k$ is not an element of $A_j \cup A'_j$, then attach an arrow to the points in $R_j$ and $R_k$ corresponding the letter, but with the arrow's head on $R_j$ and its tail on $R_k$.  Note these arrow assignments are consistent since if a single letter appears in $A_j \cup A'_j$ then it does not appear in $A_k \cup A'_k$ by the definition of word-wise partition.  Using the method just described, all arrows between distinct copies of $\R$ can be attached.

If a word $v_n = z_1z_2 \cdots z_{2M_n}$ of $p$ has a double letter, then attach arrows on $R_n$ using the following method.  As noted above, the single letters of $v_n$ in $A_n$ all have the same position parity in $z_1z_2 \cdots z_{2M_n}$ and the single letters of $v_n$ in $A'_n$ all have the other position parity by Condition (v).  Since the words of $p$ are considered up to circular permutation, assume for convenience that the single letters in $A_n$ are all in odd positions and the single letters in $A'_n$ are all in even positions.  Let $i \in E$ be a double letter of $v_n$, and suppose $z_a,z_b$ are the two occurrences of $i$ in $v_n$.  Recall that one of the indices $a,b$ is odd and the other is even by Condition (i), and that $i$ is an element of either $A_n$ or $A'_n$ by the definition of word-wise partition.  Without loss of generality, we may assume $a$ is odd and $b$ is even.  If $i \in A_n$, then attach an arrow to $R_n$ with its tail on point $a$ and its head on point $b$.  If $i \in A'_n$, then attach an arrow to $R_n$ with its head on point $a$ and its tail on point $b$.  Note that when a letter in $E$ appears in $A_n$, the tail of its corresponding arrow is on an odd point of $R_n$, and when a letter in $E$ appears in $A'_n$, the tail of its corresponding arrow is on an even point of $R_n$.  Using the method just described, all arrows with tail and head on the same copy of $\R$ can be attached.

Now to construct the virtual string $\alpha$.  The core circle of $\alpha$ corresponding to $v_n$ is the union of $R_n$ and a point at infinity with orientation extending the right-handed orientation on $R_n$.  It follows from the above methods of assigning arrows that the underlying Gauss paragraph of $\alpha$ is $p$ and the induced word-wise partition is $P$.  
\end{proof}

\begin{Lemma}\label{L:Generators}
Suppose a pair (a Gauss paragraph $p$ in an alphabet $E$, a word-wise partition $P$ of $E$) satisfies conditions (i), (iii), and (v) of the main theorem.  Let $\alpha$ be a virtual string constructed from $(p,P)$ as in Lemma \ref{L:CondVString}, so that the Gauss paragraph of $\alpha$ is $p$ and the induced word-wise partition is $P$.  Let $\rho_\alpha$ be the closed curve that is the canonical realization of $\alpha$ on the surface $\Sigma_\alpha$ as described in Section 5.  Then $H_1(\Sigma_\alpha; \Z)$ is a nontrivial free abelian group with a system of free generators consisting of three types of generators: $[S_n]$ for $1 \leq n \leq N$, $[e_i]$ for all letters $i \in E$ that are double letters of the words in $p$, and $[C_d]$ for the cyclic sequences $d \in \D_p$ associated to a subset of the single letters of the words in $p$ that will be described below.
\end{Lemma}
\begin{proof}
A graph $G_p$ can be associated to the Gauss paragraph $p$ in the following way: the words of $p$ correspond to the vertices of $G_p$, the letters in $E$ correspond to the edges of $G_p$, and the endpoints of the edge of $G_p$ corresponding to a given letter are the vertices corresponding to the words in which the letter occurs.  The graph $G_p$ is connected since by definition $p$ is not a disjoint union of other paragraphs.  Label each edge of $G_p$ with the letter in $E$ to which it corresponds.  Let $\widehat{T}$ be the set of letters that label the edges of some maximal spanning tree of $G_p$.  Notice that no double letters of the words of $p$ appear in $\widehat{T}$ since the edges labeled by double letters are loops on the vertices of $G_p$.

Although the virtual string $\alpha$ is an abstract combinatorial object, we can easily associate a trivalent directed graph $G_\alpha$ to $\alpha$.  Notice that the endpoints of the arrows of $\alpha$ divide the core circles of $\alpha$ into oriented arcs.  The edges of $G_\alpha$ correspond to these arcs and the arrows of $\alpha$, with the edges directed in ways that agree with the orientation of the arcs and a ``tail to head" orientation on the arrows.  The vertices of $G_\alpha$ correspond to the arrow endpoints of $\alpha$.  If an edge of $G_\alpha$ corresponds to an arrow of $\alpha$, label it with the letter of $E$ associated to the arrow.  The subgraph of $G_\alpha$ that consists of the edges and vertices corresponding to the core circle $S_n$ of $\alpha$ will also be referred to as $S_n$.  Note that $G_\alpha$ is not a tree since each $S_n$ is a cycle in $G_\alpha$.  Observe that the graph $G_p$ is isomorphic to the graph that results from $G_\alpha$ after contracting all $S_n$ to points and stripping directions from the remaining edges.  

Pick an edge in each $S_n$ and let $\overline{S}_n$ be the collection $S_n$ excluding the selected edge.  Let $T$ be the maximal spanning tree of $G_\alpha$ that is the union of all the $\overline{S}_n$ for $1 \leq n \leq N$ and the set of edges in $G_\alpha$ labeled with the letters in $\widehat{T}$.  The labeled edges of $G_\alpha$ not in $T$ correspond to either arrows in $\alpha$ with both endpoints on the same core circle or arrows between different core circles whose labels do not appear in $\widehat{T}$.  Therefore the labeled edges of $G_\alpha$ not in $T$ correspond to the letters in $E$ that are either double letters of some word of $p$ or single letters of words that do not appear in $\widehat{T}$, in other words the letters in $E \setminus \widehat{T}$.  There are $N$ remaining unlabeled edges of $G_\alpha$ not in $T$, namely one edge in each set $S_n \setminus \overline{S}_n$.

In general, if $X$ is a connected graph that is not a tree, then the first integral homology group of $X$ is a nontrivial free abelian group.  Given a maximal spanning tree of $X$, the first integral homology group of $X$ has a system of free generators that corresponds bijectively to the collection of edges of $X$ not in the tree.  Moreover, these generators can be calculated explicitly in a manner that we will utilize later.  These results are well known, see for example \cite{Mt}.  Recall $G_\alpha$ is a connected graph that is not a tree and $T$ is a maximal spanning tree of $G_\alpha$, so $H_1(G_\alpha) = H_1(G_\alpha;\Z)$ has a system of free generators in bijective correspondence with the collection of $\#(E \setminus \widehat{T}) + N$ edges of $G_\alpha$ not in $T$.  Since $\rho$ is a deformation retract of both $G_\alpha$ and $\Sigma = \Sigma_\alpha$, it follows that $H_1(\Sigma) = H_1(\Sigma;\Z)$ is isomorphic to $H_1(G_\alpha)$.  Therefore we can describe a system of free generators for $H_1(\Sigma)$ with $\#(E \setminus \widehat{T}) + N$ generators by investigating $H_1(G_\alpha)$.

As discussed above, there are three types of edges in $G_\alpha$ that are not in $T$.  The first type consists of unlabeled edges, with one in each set $S_n \setminus \overline{S}_n$.  Add the edge in $S_n \setminus \overline{S}_n$ to $T$ and ``prune" the resulting graph to make a cycle graph, i.e. repeatedly remove all pendant vertices until only a cycle remains.  The resulting cycle graph is $S_n$ itself.  The generator of $H_1(G_\alpha)$ represented by $S_n$ corresponds to $[S_n] \in H_1(\Sigma)$, so $[S_n]$ is in the system of free generators for $H_1(\Sigma)$.  The second type consists of edges labeled with a double letter $i \in E$ of a word $v_n$ in $p$.  Such an edge $e$ has both endpoints in $S_n$.  Add edge $e$ to $T$ and prune the resulting graph to a cycle graph.  The resulting cycle graph corresponds to the arrow $e_i$ of $\alpha$ with endpoints $a_i,b_i$ on $S_n$ and either arc $a_ib_i$ or arc $b_ia_i$ on $S_n$.  Therefore the generator of $H_1(G_\alpha)$ represented by the cycle graph corresponds to either $[e_i]$ or $[e_i]^*$ in $H_1(\Sigma)$.  We may assume the corresponding generator is $[e_i]$ since $[S_n]$ is also in the system of generators and $[e_i] + [e_i]^* = [S_n]$ in $H_1(\Sigma)$.

The third type of edges to consider are edges labeled with a single letter of two words of $p$.  Such an edge $e$ has endpoints in distinct sets $S_m$ and $S_n$.  Add edge $e$ to $T$ and prune the resulting graph to a cycle graph $C$, with orientation induced by the edge $e$.  For any of the subgraphs $S_k$ in $G_\alpha$, the intersection $C \cap S_k$ is either empty or a connected subgraph of $S_k$.  At least two such intersections must be subgraphs because $e$ has an endpoint in $S_m$ and an endpoint in $S_n$, so suppose there are $M \geq 2$ total.  Then $C$ consists of one connected subgraph in $M$ of the $S_k$ and $M$ edges that connect these subgraphs to make a cycle.  Therefore $C$ corresponds to a set of $M$ arrows between $M$ core circles of $\alpha$ and a single arc on each of these $M$ core circles.  Let $C(k)$ denote the arc of core circle $S_k$ associated to $C$ if such an arc exists, and let $C'(k)$ denote the arc in the domain of $\rho = \rho_\alpha$ corresponding to $C(k)$.  This collection of $M$ arcs is mapped to a loop $c$ on $\Sigma$ that has segments on $M$ different components of $\rho$.  The generator of $H_1(G_\alpha)$ represented by $C$ then corresponds to the equivalence class $[c] \in H_1(\Sigma)$.  The orientation of cycle graph $C$ induces an orientation on $c$.  However, on a given segment $\rho(C'(k))$ of $c$, the induced orientation may not agree with the orientation of the component $\rho(S'_k)$.  Let $\overline{c}$ be the loop on $\Sigma$ that results from removing the segment $\rho(C'(k))$ from $c$ and replacing it with the segment $\rho(S'_k \setminus (C'(k))^\circ)$.  Then $[S_k] + [c] = [\overline{c}]$ in $H_1(\Sigma)$, so we may assume the generator corresponding to $C$ is $[\overline{c}]$.  Repeating this procedure for each segment $\rho(C'(k))$ of $\bigcup_{C'(k)\ \textnormal{exists}} \rho(C'(k))$ whose induced orientation does not agree with $\rho(S'_k)$, we may conclude that the generator corresponding to $C$ is $[C_d]$ for some $d \in \D_p$.

We have now displayed a system of free generators for $H_1(\Sigma)$ with $\#(E \setminus \widehat{T}) + N$ generators, consisting of three kinds of generators: $[S_n]$ for $1 \leq n \leq N$, $[e_i]$ for all letters $i \in E$ that are double letters of words in $p$, and $[C_d]$ for sequences $d \in \D_p$ associated to the letters $i \in E \setminus \widehat{T}$ that are single letters of words in $p$.
\end{proof}
\end{section}

\begin{section}{Proof of the Main Theorem}
Now to finally prove the main theorem.

\begin{proof}[Proof of the main theorem]
Suppose the pair $(p,P)$ is realizable by a closed curve $\rho : \coprod_{n=1}^N S'_n \to S^2$, where each $S'_n$ is a copy of $S^1$.  Let $\alpha$ be the virtual string underlying $\rho$, and let $S_n$ denote the core circle of $\alpha$ corresponding to $S'_n$.  Write $p$ as a sequence $(v_1,v_2,...,v_N)$, where word $v_n$ corresponds to the core circle $S_n$.  Let $e_i = (a_i,b_i)$ for $i \in E$ be a labeling of the arrows of $\alpha$.  If a word $w$ of $p$ has a double letter $i$, then assume the $p$-sets $p_i$ and $p'_i$ contain the letters that occur exactly once in the subsequences of $w$ corresponding to the interiors of arcs $a_i b_i$ and $b_i a_i$, respectively, of the core circle of $\alpha$ associated to $w$.  Any two cycles in $S^2$ have intersection number zero, so the homological intersection form $B$ takes only the value zero in this case.  We need to show the seven conditions are satisfied.

For Condition (i), notice that two letters $i,j \in E$ are $w$-interlaced in word $w = v_m$ of $p$ if and only if the arrows $e_i,e_j$ are linked on $S_m$.  So $\#(w_i)$ is equal to the number of arrows in $arr_{m,m}(\alpha)$ that are linked with $e_i$.  Therefore $\#(w_i) \equiv n(e_i) \mtwo$.  Proposition \ref{P:BeS} implies $n(e_i) = B([e_i],[S_m]) = 0$, so $\#(w_i) \equiv n(e_i) \equiv 0 \mtwo$.  Then $\#(p_i) \equiv \sum_{k=1}^N n_{m,k}(e_i) \equiv \sum_{k \neq m} n_{m,k}(e_i) \mtwo$ and $\#(p'_i) \equiv \sum_{k=1}^N n^*_{m,k}(e_i) \equiv \sum_{k \neq m} n^*_{m,k}(e_i) \mtwo$ by Proposition \ref{P:pSetsTonij}.  But $n_{m,k}(e_i) = B([e_i],[S_k]) = 0$ and $n^*_{m,k}(e_i) = B([e_i]^*,[S_k]) = 0$ when $k \neq m$ by Proposition \ref{P:BeS}.  Hence $\#(p_i) \equiv \#(p'_{i}) \equiv 0 \mtwo$.

For Condition (ii), suppose $i \in E$ is a double letter of a word of $p$ and let $v_n$ be a different word of $p$.  Then $\#(p_i \cap o(v_n)) \equiv B([e_i],[S_n]) \equiv 0 \mtwo$ and $\#(p'_i \cap o(w)) \equiv B([e_i]^*,[S_n]) \equiv 0 \mtwo$ by Proposition \ref{P:BeS}.

For Condition (iii), suppose $p$ has at least two words.  If $v_m$ and $v_n$ are distinct words of $p$, then $B([S_m],[S_n])$ and the number of letters that these two words have in common are equivalent modulo two by Proposition \ref{P:BSS}.  Since $B([S_m],[S_n]) = 0$, it follows that the number of letters that $v_m$ and $v_n$ have in common is an even number.

For Condition (iv), if $i,j \in E$ are both double letters of the word $w$ of $p$ and are not $w$-interlaced, then Proposition \ref{P:BeeNotInt} implies $\#(w_i \cap w_j) \equiv B([e_i],[e_j]) \equiv 0 \mtwo$ since we have already shown that Condition (i) is satisfied.

For Condition (v), we must show that $P$ satisfies the two conditions in the definition of compatibility.  If $i,j \in E$ are $w$-interlaced in a word $w$ of $p$, then $w = ix_1jx_2ix_3jx_4$ for some subsequences $x_1,x_2,x_3,x_4$ of $w$.  By permuting $i$ and $j$ if necessary, we may assume the arrowtails $a_i$ and $a_j$ correspond to the first occurrences of $i$ and $j$, respectively, and the arrowheads $b_i$ and $b_j$ correspond to the second occurrences of $i$ and $j$, respectively.  If a letter $k \in p_i$, then $k$ either appears in $p'_i$ or occurs in another word of $p$.  The letter $k$ appears in $p'_i$ if and only if $k \in w_i$, so 
$$\begin{aligned}
\#(p_i) &= \#(w_i) + \#\{\textnormal{single letters of $w$ in}\ x_1jx_2\} \\
&= \#(w_i) + \#\{\textnormal{single letters of $w$ in}\ x_1\} + \#\{\textnormal{single letters of $w$ in}\ x_2\}.
\end{aligned}$$
It follows that the number of single letters in $x_1$ and the number of single letters in $x_2$ have the same parity since we have already shown that Condition (i) is satisfied.  Similar arguments for $p'_i,p_j,$ and $p'_j$ give that the numbers of single letters in $x_1,x_2,x_3,$ and $x_4$ all have the same parity.  Therefore
$\#(w_i \cap w_j) + \#\{\textnormal{single letters of $w$ in}\ x_1\} \equiv \#(w_i \cap w_j) + \#\{\textnormal{single letters of $w$ in}\ x_3\} \mtwo$.
To justify Formula 5.3.3 in the proof of Theorem 5.3.1 in \cite{Tvs}, Turaev showed that
\begin{equation}\label{E:Equation}
B([e_i],[e_j]) \equiv q(e_i,e_j) + \#\{\textnormal{single letters of}\ w\ \textnormal{in}\ x_1\} + \#(w_i \cap w_j) + 1 \mtwo
\end{equation}
when $\#(w_i) \equiv 0 \mtwo$, where $q$ is the pairing defined in Section 6.  Since $B([e_i],[e_j]) = 0$ and Condition (i) is satisfied, it follows that $$\#\{\textnormal{single letters of $w$ in}\ x_1\} + \#(w_i \cap w_j) \equiv 0 \mtwo$$ if and only if $q(e_i,e_j) \equiv 1 \mtwo$, i.e. $i$ and $j$ belong to different subsets of $P$.  Hence $P$ satisfies the first condition in the definition of compatibility.

To prove that $P$ satisfies the second condition in the definition of compatibility, suppose $i,j \in E$ are single letters in a word $w$ of $p$ that appear in the union of the two subsets associated to $w$ in $P$.  Then $w = ix_1jx_2$ for some subsequences $x_1,x_2$ of $w$.  The arrowtails of $e_i$ and $e_j$ are both on the same core circle of the virtual string, specifically the core circle corresponding to word $w$.  We have already shown that Condition (iii) is satisfied, so the length of $w$ is even by Lemma \ref{L:LengthsEven}.  Consequently $\ell(x_1) \equiv \ell(x_2) \mtwo$.  If $i$ and $j$ are in the same subset of $P$, then $q(e_i,e_j) \equiv \ell(x_1) + 1 \equiv 0 \mtwo$ and $q(e_j,e_i) \equiv \ell(x_2) + 1 \equiv 0 \mtwo$, so $\ell(x_1) \equiv \ell(x_2) \equiv 1 \mtwo$.  If $i$ and $j$ are in different subsets of $P$, then $q(e_i,e_j) \equiv \ell(x_1) + 1 \equiv 1 \mtwo$ and $q(e_j,e_i) \equiv \ell(x_2) + 1 \equiv 1 \mtwo$, so $\ell(x_1) \equiv \ell(x_2) \equiv 0 \mtwo$.  Thus $P$ is compatible with $p$.

Condition (vi) is satisfied because otherwise $B$ takes non-zero values by Proposition \ref{P:BeeDiffWords}.  Theorems \ref{T:SnDpHomology}, \ref{T:eiDpHomology}, and \ref{T:DpDpHomology} imply Condition (vii) since $B$ takes only the value zero.

Now suppose the pair $(p,P)$ satisfies all seven conditions.  Let $\alpha$ be a virtual string constructed from $(p,P)$ as in Lemma \ref{L:CondVString}, so that the Gauss word of $\alpha$ is $p$ and the induced word-wise partition is $P$.  Write $p$ as a sequence $(v_1,v_2,...,v_N)$ of words for some $N \in \mathbb{N}$ and let $S_1,S_2,...,S_N$ denote the core circles of $\alpha$, where word $v_n$ corresponds to the core circle $S_n$.  Let $e_i = (a_i,b_i)$ for $i \in E$ be a labeling of the arrows of $\alpha$.  If a word $w$ of $p$ has a double letter $i$, then assume the $p$-sets $p_i$ and $p'_i$ contain the letters that occur exactly once in the subsequences of $w$ corresponding to the interiors of arcs $a_i b_i$ and $b_i a_i$, respectively, of the core circle of $\alpha$ associated to $w$.  Using the virtual string $\alpha$, we will show that the pair $(p,P)$ is realizable by a closed curve on $S^2$.

If the surface $\Sigma = \Sigma_\alpha$ constructed in Section 5 is a 2-disc with holes, then $\Sigma$ embeds in $S^2$, so the canonical realization of $\alpha$ in $\Sigma$ gives a realization of $\alpha$ in $S^2$.  If the intersection form $B : H_1(\Sigma) \times H_1(\Sigma) \to \Z$ takes only even values, then $\Sigma$ is a disc with holes by the classification of compact surfaces.  Lemma \ref{L:Generators} implies that $H_1(\Sigma)$ is a nontrivial free abelian group with a system of free generators consisting of three types of generators: $[S_n]$ for $1 \leq n \leq N$, $[e_i]$ for all letters $i \in E$ that are double letters of words of $p$, and $[C_d]$ for the cyclic sequences $d \in \D_p$ associated to a subset of the single letters of the words in $p$ as described in the proof of the lemma.  To show that $B$ takes only even values on $H_1(\Sigma)$, it suffices to show $B$ takes only even values on these generators.

For three cases, Turaev's arguments in the proof of Theorem 5.3.1 in \cite{Tvs} can be generalized as follows.  If $i \in E$ is a double letter of the word $w = v_m$ of $p$, then $B([e_i],[S_m]) = n(e_i)$ by Proposition \ref{P:BeS}.  But $n(e_i) \equiv \#(w_i) \mtwo$, so $B([e_i],[S_m])$ is even by Condition (i).  Suppose $j \in E$ is another double letter of word $w$.  If the arrows $e_i$ and $e_j$ are not linked, then $B([e_i],[e_j]) \equiv \#(w_i \cap w_j) \equiv 0 \mtwo$ by Proposition \ref{P:BeeNotInt} and conditions (i) and (iv).  If $e_i$ and $e_j$ are linked, then first assume that their endpoints lie in the cyclic order $a_i,a_j,b_i,b_j$ around $S_m$, with the subsequence $x$ of $w$ corresponding to the interior of the arc $a_ia_j$ on $S_m$.  If $i,j$ appear in the same subset of $P$, then $q(e_i,e_j) \equiv 0 \mtwo$ by the construction of $\alpha$ and $\#(w_i \cap w_j) + \#\{\textnormal{single letters in}\ x\} \equiv 1 \mtwo$ by Condition (v).  Therefore $B([e_i],[e_j])$ is even by Formula (\ref{E:Equation}) and Condition (i).  If $i,j$ appear in different subsets of $P$, then $q(e_i,e_j) \equiv 1 \mtwo$ by the construction of $\alpha$ and $\#(w_i \cap w_j) + \#\{\textnormal{single letters in}\ x\} \equiv 0 \mtwo$ by Condition (v).  Therefore $B([e_i],[e_j])$ is even by Formula (\ref{E:Equation}) and Condition (i).  The case where the arrows $e_i$ and $e_j$ are linked and their endpoints are assumed to lie in the cyclic order $a_i,b_j,b_i,a_j$ around $S_m$ follows from the previous case and the skew-symmetry of $B$. 

We will now consider all other cases.  If $i \in E$ is a double letter of the word $v_m$ of $p$ and $j \in E$ is a double letter of a different word $v_n$ of $p$, then $B([e_i],[e_j]) \equiv \#(p_i \cap p_j) \equiv 0 \mtwo$ by Proposition \ref{P:BeeDiffWords} and Condition (vi).  Moreover, $B([e_i],[S_n]) \equiv \#(p_i \cap o(v_n)) \equiv 0 \mtwo$ by Proposition \ref{P:BeS} and Condition (ii).  Proposition \ref{P:BSS} implies $B([S_m],[S_n])$ is equivalent modulo two to the number of letters that the words $v_m$ and $v_n$ of $p$ have in common, so $B([S_m],[S_n]) \equiv 0 \mtwo$ by Condition (iii).  Since $P$ is compatible with $\D_p$ by Condition (vii), it follows from theorems \ref{T:SnDpHomology}, \ref{T:eiDpHomology}, and \ref{T:DpDpHomology} that
$$B([S_n],[C_{d_1}]) \equiv B([e_i],[C_{d_1}]) \equiv B([C_{d_1}],[C_{d_2}]) \equiv 0 \mtwo$$
for all $1 \leq n \leq N$, double letters $i \in E$ of the words in $p$, and $d_1,d_2 \in \D_p$.

We have proved that $B$ takes only even values on $H_1(\Sigma)$, which implies $\Sigma$ is a disc with holes.  Then $\Sigma$ embeds in $S^2$, giving a realization of $\alpha$ in $S^2$ induced by the canonical realization of $\alpha$ in $\Sigma$.  Thus the pair $(p,P)$ is realizable by a closed curve on $S^2$.  
\end{proof}

\begin{remark}
It may be interesting to write a computer program that implements an algorithm for checking the seven conditions of the main theorem.  We already have some specific ideas about such an algorithm, but these ideas will be postponed until a future paper.
\end{remark}
\end{section}

\begin{section}{Additional Results}
Throughout this section, let $p$ be a Gauss paragraph in an alphabet $E$.
\begin{Lemma}\label{L:ppwwout}
Suppose $i,j \in E$ are $w$-interlaced in a word $w$ of $p$.  Write $w$ in the form $i\ x_1\ j\ x_2\ i\ x_3\ j\ x_4$ for some subsequences $x_1,x_2,x_3,x_4$ of $w$, and assume $p_i = o(i\ x_1\ j\ x_2\ i)$, $p'_i = o(i\ x_3\ j\ x_4\ i)$, $p_j = o(j\ x_2\ i\ x_3\ j)$, and $p'_j = o(j\ x_4\ i\ x_1\ j)$.  Then
$$\begin{aligned}
\#(p_i \cap p_j) &= \#(w_i \cap w_j) + out(x_2) \\
\#(p_i \cap p'_j) &= \#(w_i \cap w_j) + out(x_1) \\
\#(p'_i \cap p_j) &= \#(w_i \cap w_j) + out(x_3) \\
\#(p'_i \cap p'_j) &= \#(w_i \cap w_j) + out(x_4),
\end{aligned}$$
where $out(x_k) \in \N$ denotes the number of letters that the subsequence $x_k$ has in common with the other words of $p$.
\end{Lemma}
\begin{proof}
Note $\#(p_i \cap p_j)$ is the number of letters that occur exactly once in both $p_i$ and $p_j$.  There are three types of such letters: letters in both $x_1$ and $x_3$, letters in both $x_2$ and $x_4$, and letters in $out(x_2)$.  A letter appears in both $x_1$ and $x_3$ or both $x_2$ and $x_4$ if and only if it appears in $\#(w_i \cap w_j)$.  Hence $\#(p_i \cap p_j) = \#(w_i \cap w_j) + out(x_2)$.  The remaining three claims can be proved similarly.
\end{proof}

Note that when $p$ is a one-word Gauss paragraph, the above lemma implies
$\#(w_i \cap w_j) = \#(p_i \cap p_j) = \#(p_i \cap p'_j) = \#(p'_i \cap p_j) = \#(p'_i \cap p'_j)$
since $out(x_2) = out(x_1) = out(x_3) = out(x_4) = 0$.

\begin{Lemma}\label{L:pwout}
Suppose $i \in E$ occurs twice in a word $w$ of $p$.  Write $w$ in the form $i x_1 i x_2$ for some subsequences $x_1,x_2$ of $w$, and assume $p_i = o(i x_1 i)$ and $p'_i = o(i x_2 i)$.  Then
$$\begin{aligned}
\#(p_i) &= \#(w_i) + out(x_1) \\
\#(p'_i) &= \#(w_i) + out(x_2)
\end{aligned}$$
\end{Lemma}
\begin{proof}
If $j \in p_i$, then $j$ either appears in $p'_i$ or occurs in another word of $p$.  The letter $j$ appears in $p'_i$ if and only if $j \in w_i$.  The letter $j$ occurs in another word of $p$ if and only if $j \in out(x_1)$.  The equation involving $p'_i$ follows from similar arguments.
\end{proof}

\begin{Thm}\label{T:Equivalences}
Suppose $i,j \in E$ are $w$-interlaced in a word $w$ of $p$.  Write $w$ in the form $i\ x_1\ j\ x_2\ i\ x_3\ j\ x_4$ for some subsequences $x_1,x_2,x_3,x_4$ of $w$, and assume $p_i = o(i\ x_1\ j\ x_2\ i)$ and $p'_i = o(i\ x_3\ j\ x_4\ i)$.  Then the following are equivalent:
\begin{enumeratei}
\item{$\#(p_i \cap p_j) \equiv \#(p_i \cap p'_j) \equiv \#(p'_i \cap p_j) \equiv \#(p'_i \cap p'_j) \mtwo;$}
\item{$out(x_1) \equiv out(x_2) \equiv out(x_3) \equiv out(x_4) \mtwo;$}
\item{$out(x_1 j x_2) \equiv out(x_2 i x_3) \equiv out(x_3 j x_4) \equiv out(x_4 i x_1) \equiv 0 \mtwo;$}
\item{$\#(p_i) \equiv \#(p'_i) \equiv \#(w_i) \mtwo$ and $\#(p_j) \equiv \#(p'_j) \equiv \#(w_j) \mtwo.$}
\end{enumeratei}
\end{Thm}
\begin{proof}
Statements (i) and (ii) are equivalent via Lemma \ref{L:ppwwout}.  Statements (ii) and (iii) are equivalent since $i$ and $j$ are double letters of the word $w$.  Statements (iii) and (iv) are equivalent via Lemma \ref{L:pwout}.
\end{proof}

\begin{Prop}\label{P:pEvens}
Suppose $i \in E$ is a double letter of a word $w$ in $p$ and $j \in E$ is a double letter of a different word $w'$ in $p$.  If
$$\#(p_i \cap p_j) \equiv \#(p_i \cap p'_j) \equiv \#(p'_i \cap p_j) \equiv \#(p'_i \cap p'_j) \mtwo,$$
then $p_i$ and $p'_i$ each have an even number of letters in common with word $w'$, and hence the words $w$ and $w'$ have an even number of letters in common.
\end{Prop}
\begin{proof}
Note $\#(p_i \cap p_j) + \#(p_i \cap p'_j) \equiv \#(p'_i \cap p_j) + \#(p'_i \cap p'_j) \equiv 0 \mtwo$ by hypothesis.  The first and second sums are equal to the number of letters that $p_i$ and $p'_i$, respectively, have in common with word $w'$.  If both $p_i$ and $p'_i$ have an even number of letters in common with $w'$, then $w$ must have an even number of letters in common with $w'$.
\end{proof}
\end{section}

\begin{section}{Acknowledgements}
I wish to thank my advisor R.A. Litherland for the many valuable conversations and suggestions, and especially for his patience.  
\end{section}

\bibliographystyle{abbrv}
\bibliography{references}

\begin{thebibliography}{1}

\bibitem{Ccic}
J.~S. Carter.
\newblock Classifying immersed curves.
\newblock {\em Proc. Amer. Math. Soc.}, 111(1):281--287, 1991.

\bibitem{DTcok}
C.~H. Dowker and M.~B. Thistlethwaite.
\newblock Classification of knot projections.
\newblock {\em Topology Appl.}, 16(1):19--31, 1983.

\bibitem{Gw}
C.~F. Gauss.
\newblock {\em Werke}, pages 272 and 282--286.
\newblock Teubner, Leipzig, 1900.

\bibitem{LMafs}
L.~Lov{\'a}sz and M.~L. Marx.
\newblock A forbidden substructure characterization of {G}auss codes.
\newblock {\em Acta Sci. Math. (Szeged)}, 38(1--2):115--119, 1976.

\bibitem{Mt}
J.~R. Munkres.
\newblock {\em Topology: a first course}.
\newblock Prentice-Hall Inc., Englewood Cliffs, N.J., second edition, 2000.

\bibitem{Rsad}
P.~Rosenstiehl.
\newblock Solution alg\'ebrique du probl\`eme de {G}auss sur la permutation des
  points d'intersection d'une ou plusieurs courbes ferm\'ees du plan.
\newblock {\em C. R. Acad. Sci. Paris S\'er. A-B}, 283(8):Ai, A551--A553, 1976.

\bibitem{Tvs}
V.~Turaev.
\newblock Virtual strings.
\newblock {\em ArXiV:math.GT/0310218}, 2003.
\newblock preprint.

\end{thebibliography}

\end{document}